\documentclass{article}
\usepackage{graphicx}
\usepackage{epstopdf}
\usepackage{indentfirst}

\usepackage[latin1]{inputenc}
\usepackage{tikz}
\usetikzlibrary{trees}
\usepackage{verbatim}

\tikzstyle{level 1}=[level distance=3cm, sibling distance=3cm]
\tikzstyle{level 2}=[level distance=2cm, sibling distance=2cm]

\tikzstyle{probranch} = [circle, minimum width=5pt,fill=black, inner sep=0pt]

\makeatletter
\g@addto@macro\Gin@extensions{,.gif}
\makeatother
\begin{document}

\title{Singular Points of High Multiplicity for Septic Curves
}
\author{Nicholas J. Willis, David A. Weinberg}
\maketitle

\abstract{For real irreducible algebraic curves of the seventh degree, there are 22 types of singular points of multiplicity six, 174 types of singular points of multiplicity five, and at least 182 types of singular points of multiplicity four. For complex irreducible algebraic curves of the seventh degree, there are 12 types of singular points of multiplicity six, 92 types of singular points of multiplicity five, and at least 92 types of singular points of multiplicity four. In the final section of the paper, a wide variety of open problems on the classification of singular points of plane algebraic cuves is explicitly described. }
\section{Introduction}

In this paper singular points of multiplicities 4, 5, and 6 for real irreducible seventh degree curves are 
studied.  In previous papers, the authors have classified singular points for curves of degrees 4, 5, and 6 [3],[4],[5]. \\

Here is a table indicating the number of types of singular points for various families of curves: \\ \\

\begin{tabular}{l|c}
\hline
Family of curves & Number of types of singular points\\ \hline  \\
Irreducible real quartics & 13 \\
Reducible real quartics & 17 \\
Irreducible complex quartics & 9 \\
Reducible complex quartics & 10 \\
Irreducible real quintics & 42 \\
Reducible real quintics & 49 \\
Irreducible complex quintics & 28 \\
Reducible complex quintics & 33 \\
Irreducible real sextics & 191 \\
Reducible real sextics & 189? \\
Irreducible complex sextics & 108 \\
Reducible complex sextics & 106 \\
\hline
\end{tabular}
\\ \\
\\

\textbf{Notational convention.} 
The reader may have noticed a question mark in the above table under the entry for reducible real sextics. This means that 189 types were constructed, but we can neither prove that there are no more types, nor can we construct further types due to the current limitations (as of 2014) on memory in desktop computing. This same usage of notation will occur later in this paper for several types of tree diagrams for points of multiplicity four. \\

To the best of our knowledge, there has been no systematic study of singular points of seventh degree curves until now.  The definition of the equivalence relation we study for singular points has been described and explained in detail in previous papers.  See [3] or [4].  We follow the same conventions that were described in those papers. (The reader is warned that because of some clerical error, the 2nd diagram on p. 102 and the diagram at the bottom of p. 103 are incorrect in [3].) \\ 

\underline{Equivalence relation for complex curves} \\

Two singular points of complex algebraic curves are equivalent if their pro-branches (Puiseux expansions) have the same exponents of contact. \\ 

In [3] and [4] we described in detail how to assign a weighted tree diagram to each singular point.  Such a diagram nicely exhibits all the information about the exponents of contact.  Pro-branches and exponents of contact were defined in [2], p.68.  Let us further clarify here that if the equation of the curve has real coefficients, then it is possible to have two or more expansions whose coefficients are complex conjugates in pairs of expansions and such that the only real point corresponding to the expansions is the origin.  (The simplest possible example is $y^2+x^2=0$ : the Puiseux expansions are $y=\pm ix$, but the only real point is $(0,0)$.)  In such cases, we refer to the exponents of contact as *-exponents of contact and indicate this in the tree diagram by braces.  Thus, \\ 

\underline{Equivalence relation for real curves} \\ 

Two singular points of real algebraic curves are equivalent if their pro-branches have the same exponents of contact and the same *-exponents of contact. \\

The classification proceeds one multiplicity at a time.  For a given multiplicity, we first choose the tangent cone.  For each tangent cone, we consider all possible Newton polygons.  This then yields several families of curves.  For some of these families, a sequence of Puiseux expansions is computed using computer algebra in the form of Maple.  After each Puiseux expansion, the family is modified by a polynomial condition obtained during the Puiseux expansion which is obtained by computer algebra.  The sequence of calculations stops when the family consists entirely of reducible curves. (Reducible curves are treated in a separate study.) \\ 

In the next section, we exhibit, for each multiplicity, the possible tangent cones, the Newton polygons, and the associated families of curves.  For each family of curves, we indicate all possible types of singular points by exhibiting the diagrams.  As we explained in previous papers [3] and [4], Puiseux expansion calculations are only necessary when there is a multiple root along any segment of the Newton polygon. Numerous detailed examples of such computer calculations are shown in [3] and [4]. Thus, in what follows, the reader is expected to know when computer calculations have been performed and when the diagrams were obtained directly from the Newton polygon without computer calculation. Since we have presented an exhaustive list of families below, the interested reader can easily check any case using, for example, Maple computer algebra. \\

\textbf{Convention about Newton polygons and nonzero coefficients.} 
When a Newton polygon is exhibited, it will be assumed that the associated family has nonzero coefficients corresponding to the endpoints of the segments on the Newton polygon. For example, for the very first Newton polygon and family in the next section, it is assumed that $a\neq0$ without explicit mention.  One notable convention is that Newton-Polygons will be in the style of Walker[1]\\

Another feature of our presentation that needs to be explained is that the diagrams for the singular points are only listed once. For example, if some of the diagrams yielded by a given family were already listed in a previous case, then only the new diagrams for that family will be listed. A more interesting instance of this is that in some special cases, by an expedient change of variables, the Newton polygon can be transformed, reducing the study of one family to another, and thus saving a large amount of computation. \\

In a couple of cases, this technique can be used to eliminate diagram types because the transformed Newton polygon does not exist for seventh degree curves. The interested reader is expected to identify and verify the above situations. \\

\textbf{Acknowledgment}  The authors wish to thank Esther Hunt for many valuable questions and comments and Janelle Dockter for assistance with Latex.

\section{Real Singular Points of Multiplicity 6}

\noindent Tangent Cone: $y^6$ \\ 
\\
Newton Polygon 1:\\  
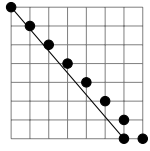
\begin{center}

\end{center}
\noindent
Family: \\ $y^6+ax^7+bx^6y+cx^5y^2+dx^4y^3+ex^3y^4+fx^2y^5+gxy^6+hy^7=0$ \\ 
\noindent \\
Diagram: \\
\begin{center}
\tikzstyle{level 1}=[level distance=9mm, sibling distance=3mm]
\tikzstyle{level 2}=[level distance=6mm, sibling distance=6mm]
\begin{tabular}{l r}
\begin{tikzpicture}[grow=right, sloped]
\draw  [-, rounded corners, thick] (1,-.8) -- (1.2,-.6) -- (1,-.4);
\draw  [-, rounded corners, thick] (1,-.15) -- (1.2,0) -- (1,.15);
\node[probranch] {}
    child {
        node[probranch] {}        
        edge from parent 
            node[] {}
            node[]  {}
    }
    child {
        node[probranch] {}        
        edge from parent 
            node[] {}
            node[]  {}
    }
    child {
        node[probranch] {}        
        edge from parent 
            node[] {}
            node[]  {}
    }
    child {
        node[probranch] {}        
        edge from parent 
            node[] {}
            node[]  {}
    }
    child {
        node[probranch] {}        
        edge from parent 
            node[] {}
            node[]  {}
    }
    child {
        node[probranch, label=above: $7/6$] {}        
        edge from parent         
            node[] {}
            node[]  {}
    };
\end{tikzpicture}\\
\end{tabular}%
\end{center}
Family: \\ $y^5(y-x)+ax^7+bx^6y+cx^5y^2+dx^4y^3+ex^3y^4+fx^2y^5+gxy^6+hy^7=0$ \\ 
\\
Diagram: \\
\begin{center}
\tikzstyle{level 1}=[level distance=9mm, sibling distance=9mm]
\tikzstyle{level 2}=[level distance=6mm, sibling distance=3mm]
\begin{tabular}{l r}
\begin{tikzpicture}[grow=right, sloped]
\draw  [-, rounded corners, thick] (1.6,1.1) -- (1.8,.9) -- (1.6,.7);
\draw  [-, rounded corners, thick] (1.6,.5) -- (1.8,.3) -- (1.6,.1);
\node[probranch] {}
    child {
        node[probranch] {}
        		child {
        				node[probranch] {}    
        				edge from parent 
           					node[] {}
            				node[]  {}
    				}    
        edge from parent 
            node[] {}
            node[]  {}
    }
    child {
        node[probranch, label=above: $1$] {}
        		child {
        				node[probranch] {}    
        				edge from parent 
           					node[] {}
            				node[]  {}
    				}
    				child {
        				node[probranch] {}    
        				edge from parent 
           					node[] {}
            				node[]  {}
    				}
    				child {
        				node[probranch] {}    
        				edge from parent 
           					node[] {}
            				node[]  {}
    				}
    				child {
        				node[probranch] {}    
        				edge from parent 
           					node[] {}
            				node[]  {}
    				}      
    				child {
        				node[probranch, label=above: $6/5$] {}    
        				edge from parent 
           					node[] {}
            				node[]  {}
    				}          
        edge from parent         
            node[] {}
            node[]  {}
    };
\end{tikzpicture} \\
\end{tabular}%
\end{center}
Families: \\ $y^4(y-x)^2+ax^7+bx^6y+cx^5y^2+dx^4y^3+ex^3y^4+fx^2y^5+gxy^6+hy^7=0$, \\ \\ $y^4(y-x)(y-2x)+ax^7+bx^6y+cx^5y^2+dx^4y^3+ex^3y^4+fx^2y^5+gxy^6+hy^7=0$, \\ \\ $y^4(y^2+x^2)+ax^7+bx^6y+cx^5y^2+dx^4y^3+ex^3y^4+fx^2y^5+gxy^6+hy^7=0$ \\ 
\\
Diagrams: \\
\begin{center}
\tikzstyle{level 1}=[level distance=9mm, sibling distance=6.5mm]
\tikzstyle{level 2}=[level distance=6mm, sibling distance=2.5mm]
\begin{tabular}{l r}
\begin{tikzpicture}[grow=right, sloped]
\draw  [-, rounded corners, thick] (1.6,1.05) -- (1.8,.9) -- (1.6,.75);
\node[probranch] {}
    child {
        node[probranch] {}
        		child {
        				node[probranch] {}    
        				edge from parent 
           					node[] {}
            				node[]  {}
    				}    
        edge from parent 
            node[] {}
            node[]  {}
    }
    child {
        node[probranch] {}
        		child {
        				node[probranch] {}    
        				edge from parent 
           					node[] {}
            				node[]  {}
    				}    
        edge from parent 
            node[] {}
            node[]  {}
    }
    child {
        node[probranch, label=above: $1$] {}
        		child {
        				node[probranch] {}    
        				edge from parent 
           					node[] {}
            				node[]  {}
    				}
    				child {
        				node[probranch] {}    
        				edge from parent 
           					node[] {}
            				node[]  {}
    				}
    				child {
        				node[probranch] {}    
        				edge from parent 
           					node[] {}
            				node[]  {}
    				}      
    				child {
        				node[probranch, label=above: $5/4$] {}    
        				edge from parent 
           					node[] {}
            				node[]  {}
    				}          
        edge from parent         
            node[] {}
            node[]  {}
    };
\end{tikzpicture}\\
\end{tabular}
\tikzstyle{level 1}=[level distance=9mm, sibling distance=6.5mm]
\tikzstyle{level 2}=[level distance=6mm, sibling distance=2.5mm]
\begin{tabular}{l r}
\begin{tikzpicture}[grow=right, sloped]
\draw  [-, rounded corners, thick] (1.6,1.05) -- (1.8,.9) -- (1.6,.75);
\draw [-, rounded corners, thick] (1,-.05) -- (1.2,-.3) -- (1,-.6);
\node[probranch] {}
    child {
        node[probranch] {}
        		child {
        				node[probranch] {}    
        				edge from parent 
           					node[] {}
            				node[]  {}
    				}    
        edge from parent 
            node[] {}
            node[]  {}
    }
    child {
        node[probranch] {}
        		child {
        				node[probranch] {}    
        				edge from parent 
           					node[] {}
            				node[]  {}
    				}    
        edge from parent 
            node[] {}
            node[]  {}
    }
    child {
        node[probranch, label=above: $1$] {}
        		child {
        				node[probranch] {}    
        				edge from parent 
           					node[] {}
            				node[]  {}
    				}
    				child {
        				node[probranch] {}    
        				edge from parent 
           					node[] {}
            				node[]  {}
    				}
    				child {
        				node[probranch] {}    
        				edge from parent 
           					node[] {}
            				node[]  {}
    				}      
    				child {
        				node[probranch, label=above: $5/4$] {}    
        				edge from parent 
           					node[] {}
            				node[]  {}
    				}          
        edge from parent         
            node[] {}
            node[]  {}
    };
\end{tikzpicture} \\
\end{tabular}
\tikzstyle{level 1}=[level distance=9mm, sibling distance=9mm]
\tikzstyle{level 2}=[level distance = 6mm, sibling distance=3mm]
\tikzstyle{level 3}=[level distance=6mm, sibling distance=3mm]
\begin{tabular}{l r}
\begin{tikzpicture}[grow=right, sloped]
\draw  [-, rounded corners, thick] (2.2,.95) -- (2.4,.75) -- (2.2,.55);
\node[probranch] {}
    child {
        node[probranch] {}
        		child {
        				node[probranch] {}
        						child {
        								node[probranch] {}    
        								edge from parent 
           									node[] {}
            								node[]  {}
    								} 
    								child {
        								node[probranch] {}    
        								edge from parent 
           									node[] {}
            								node[]  {}
    								}     
        				edge from parent 
           					node[] {}
            				node[]  {}
    				}    
        edge from parent 
            node[] {}
            node[]  {}
    }
    child {
        node[probranch, label=above: $1$] {}
        		child {
        				node[probranch] {}
        						child {
        								node[probranch] {}    
        								edge from parent 
           									node[] {}
            								node[]  {}
    								}     
        				edge from parent 
           					node[] {}
            				node[]  {}
    				}
    				child {
        				node[probranch] {}
        						child {
        								node[probranch] {}    
        								edge from parent 
           									node[] {}
            								node[]  {}
    								}     
        				edge from parent 
           					node[] {}
            				node[]  {}
    				}
    				child {
        				node[probranch] {}
        						child {
        								node[probranch] {}    
        								edge from parent 
           									node[] {}
            								node[]  {}
    								}     
        				edge from parent 
           					node[] {}
            				node[]  {}
    				}      
    				child {
        				node[probranch, label=above: $5/4$] {}
        						child {
        								node[probranch, label=above: $3/2$] {}    
        								edge from parent 
           									node[] {}
            								node[]  {}
    								}     
        				edge from parent 
           					node[] {}
            				node[]  {}
    				}          
        edge from parent         
            node[] {}
            node[]  {}
    };
\end{tikzpicture} \\
\end{tabular}%
\end{center}
Families: \\ $y^3(y-x)^3+ax^7+bx^6y+cx^5y^2+dx^4y^3+ex^3y^4+fx^2y^5+gxy^6+hy^7=0$, \\ \\ $y^3(y-x)^2(y-2x)+ax^7+bx^6y+cx^5y^2+dx^4y^3+ex^3y^4+fx^2y^5+gxy^6+hy^7=0$, \\ \\ $y^3(y-x)(y-2x)(y-ax)+bx^7+cx^6y+dx^5y^2+ex^4y^3+fx^3y^4+gx^2y^5+hxy^6+jy^7=0$ where $a\neq0,1,2$, \\ \\ $y^3(y-x)(y^2+x^2)+ax^7+bx^6y+cx^5y^2+dx^4y^3+ex^3y^4+fx^2y^5+gxy^6+hy^7=0$ \\ 
\\
Diagrams: 
\begin{center}
\tikzstyle{level 1}=[level distance=9mm, sibling distance=6mm]
\tikzstyle{level 2}=[level distance=6mm, sibling distance=3mm]
\begin{tabular}{l r}
\begin{tikzpicture}[grow=right, sloped]
\draw  [-, rounded corners, thick] (1.6,.9) -- (1.8,.75) -- (1.6,.6);
\node[probranch] {}
    child {
        node[probranch] {}
        		child {
        				node[probranch] {}    
        				edge from parent 
           					node[] {}
            				node[]  {}
    				}    
        edge from parent 
            node[] {}
            node[]  {}
    }
    child {
        node[probranch] {}
        		child {
        				node[probranch] {}    
        				edge from parent 
           					node[] {}
            				node[]  {}
    				}    
        edge from parent 
            node[] {}
            node[]  {}
    }
    child {
        node[probranch] {}
        		child {
        				node[probranch] {}    
        				edge from parent 
           					node[] {}
            				node[]  {}
    				}    
        edge from parent 
            node[] {}
            node[]  {}
    }
    child {
        node[probranch, label=above: $1$] {}
        		child {
        				node[probranch] {}    
        				edge from parent 
           					node[] {}
            				node[]  {}
    				}
    				child {
        				node[probranch] {}    
        				edge from parent 
           					node[] {}
            				node[]  {}
    				}     
    				child {
        				node[probranch, label=above: $4/3$] {}    
        				edge from parent 
           					node[] {}
            				node[]  {}
    				}          
        edge from parent         
            node[] {}
            node[]  {}
    };
\end{tikzpicture}\\
\end{tabular}
\tikzstyle{level 1}=[level distance=9mm, sibling distance=6mm]
\tikzstyle{level 2}=[level distance=6mm, sibling distance=3mm]
\begin{tabular}{l r}
\begin{tikzpicture}[grow=right, sloped]
\draw  [-, rounded corners, thick] (1.6,.9) -- (1.8,.75) -- (1.6,.6);
\draw [-, rounded corners, thick] (1, .2) -- (1.2, 0) -- (1, -.2);
\node[probranch] {}
    child {
        node[probranch] {}
        		child {
        				node[probranch] {}    
        				edge from parent 
           					node[] {}
            				node[]  {}
    				}    
        edge from parent 
            node[] {}
            node[]  {}
    }
    child {
        node[probranch] {}
        		child {
        				node[probranch] {}    
        				edge from parent 
           					node[] {}
            				node[]  {}
    				}    
        edge from parent 
            node[] {}
            node[]  {}
    }
    child {
        node[probranch] {}
        		child {
        				node[probranch] {}    
        				edge from parent 
           					node[] {}
            				node[]  {}
    				}    
        edge from parent 
            node[] {}
            node[]  {}
    }
    child {
        node[probranch, label=above: $1$] {}
        		child {
        				node[probranch] {}    
        				edge from parent 
           					node[] {}
            				node[]  {}
    				}
    				child {
        				node[probranch] {}    
        				edge from parent 
           					node[] {}
            				node[]  {}
    				}     
    				child {
        				node[probranch, label=above: $4/3$] {}    
        				edge from parent 
           					node[] {}
            				node[]  {}
    				}          
        edge from parent         
            node[] {}
            node[]  {}
    };
\end{tikzpicture} \\
\end{tabular}
\tikzstyle{level 1}=[level distance=9mm, sibling distance=9mm]
\tikzstyle{level 2}=[level distance=6mm, sibling distance=3mm]
\begin{tabular}{l r}
\begin{tikzpicture}[grow=right, sloped]
\draw  [-, rounded corners, thick] (1.6,.15) -- (1.8,.3) -- (1.6,.45);
\draw [-, rounded corners, thick] (1.6, -.45) -- (1.8,-.6) -- (1.6,-.75);
\node[probranch] {}
    child {
        node[probranch] {}
        		child {
        				node[probranch] {}    
        				edge from parent 
           					node[] {}
            				node[]  {}
    				}
    				child {
        				node[probranch] {}    
        				edge from parent 
           					node[] {}
            				node[]  {}
    				}
    				child {
        				node[probranch] {}    
        				edge from parent 
           					node[] {}
            				node[]  {}
    				}    
        edge from parent 
            node[] {}
            node[]  {}
    }
    child {
        node[probranch, label=above: $1$] {}
        		child {
        				node[probranch] {}    
        				edge from parent 
           					node[] {}
            				node[]  {}
    				}
    				child {
        				node[probranch] {}    
        				edge from parent 
           					node[] {}
            				node[]  {}
    				}      
    				child {
        				node[probranch, label=above: $4/3$] {}    
        				edge from parent 
           					node[] {}
            				node[]  {}
    				}          
        edge from parent         
            node[] {}
            node[]  {}
    };
\end{tikzpicture} \\
\end{tabular}
\tikzstyle{level 1}=[level distance=9mm, sibling distance=9mm]
\tikzstyle{level 2}=[level distance=6mm, sibling distance=3mm]
\tikzstyle{level 3}=[level distance=6mm, sibling distance=3mm]
\begin{tabular}{l r}
\begin{tikzpicture}[grow=right, sloped]
\draw  [-, rounded corners, thick] (1.6,.5) -- (1.8,.6) -- (1.6,.7);
\node[probranch] {}
    child {
        node[probranch] {}
        		child {
        				node[probranch] {}
        						child {
        								node[probranch] {}    
        								edge from parent 
           									node[] {}
            								node[]  {}
    								} 
    								child {
        								node[probranch] {}    
        								edge from parent 
           									node[] {}
            								node[]  {}
    								}     
								child {
        				node[probranch] {}
        				edge from parent 
           					node[] {}
            				node[]  {}
    				}
        				edge from parent 
           					node[] {}
            				node[]  {}
    				}    
        edge from parent 
            node[] {}
            node[]  {}
    }
    child {
        node[probranch, label=above: $1$] {}
    				child {
        				node[probranch] {}
        						child {
        								node[probranch] {}    
        								edge from parent 
           									node[] {}
            								node[]  {}
    								}     
        				edge from parent 
           					node[] {}
            				node[]  {}
    				}
    				child {
        				node[probranch] {}
        						child {
        								node[probranch] {}    
        								edge from parent 
           									node[] {}
            								node[]  {}
    								}     
        				edge from parent 
           					node[] {}
            				node[]  {}
    				}      
    				child {
        				node[probranch, label=above: $4/3$] {}
        						child {
        								node[probranch, label=above: $3/2$] {}    
        								edge from parent 
           									node[] {}
            								node[]  {}
    								}     
        				edge from parent 
           					node[] {}
            				node[]  {}
    				}          
        edge from parent         
            node[] {}
            node[]  {}
    };
\end{tikzpicture} \\
\end{tabular}%
\end{center}
Families: \\ $y^2(y-x)(y+x)(y+ax)(y+bx)+cx^7+dx^6y+ex^5y^2+fx^4y^3+gx^3y^4+hx^2y^5+jxy^6+ky^7=0$, \\ \\ $y^2(y-x)^2(y+x)(y+ax)+bx^7+cx^6y+dx^5y^2+ex^4y^3+fx^3y^4+gx^2y^5+hxy^6+jy^7=0$ where $a\neq0,1,-1$, \\ \\ $y^2(y-x)^2(y+x)^2+ax^7+bx^6y+cx^5y^2+dx^4y^3+ex^3y^4+fx^2y^5+gxy^6+hy^7=0$, \\ \\ $y^2(y^2+x^2)(y+ax)(y+bx)+cx^7+dx^6y+ex^5y^2+fx^4y^3+gx^3y^4+hx^2y^5+jxy^6+ky^7=0$, where $a, b\neq0$ and $a\neq b$\\  \\$y^2(y^2+x^2)(y+ax)^2+bx^7+cx^6y+dx^5y^2+ex^4y^3+fx^3y^4+gx^2y^5+hxy^6+jy^7=0$ where $a\neq0$,\\ \\ $y^2(y^2+x^2)(y^2+ax^2)+ax^7+bx^6y+cx^5y^2+dx^4y^3+ex^3y^4+fx^2y^5+gxy^6+hy^7=0$ where $a>0$ and $a\neq1$, \\ \\ $y^2(y^2+x^2)^2+ax^7+bx^6y+cx^5y^2+dx^4y^3+ex^3y^4+fx^2y^5+gxy^6+hy^7=0$, \\ 
\\
Diagrams: 
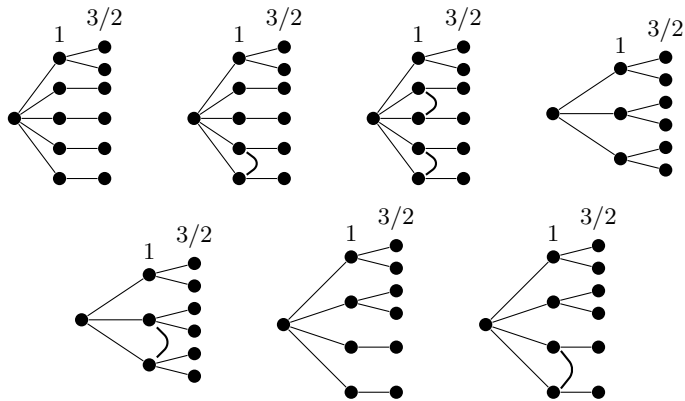
\begin{center}
\tikzstyle{level 1}=[level distance=6mm, sibling distance=4mm]
\tikzstyle{level 2}=[level distance=6mm, sibling distance=3mm]
\begin{tabular}{l r}
\begin{tikzpicture}[grow=right, sloped]
\node[probranch] {}
    child {
        node[probranch] {}
        		child {
        				node[probranch] {}    
        				edge from parent 
           					node[] {}
            				node[]  {}
    				}    
        edge from parent 
            node[] {}
            node[]  {}
    }
    child {
        node[probranch] {}
        		child {
        				node[probranch] {}    
        				edge from parent 
           					node[] {}
            				node[]  {}
    				}    
        edge from parent 
            node[] {}
            node[]  {}
    }
    child {
        node[probranch] {}
        		child {
        				node[probranch] {}    
        				edge from parent 
           					node[] {}
            				node[]  {}
    				}    
        edge from parent 
            node[] {}
            node[]  {}
    }
    child {
        node[probranch] {}
        		child {
        				node[probranch] {}    
        				edge from parent 
           					node[] {}
            				node[]  {}
    				}    
        edge from parent 
            node[] {}
            node[]  {}
    }
    child {
        node[probranch, label=above: $1$] {}
        		child {
        				node[probranch] {}    
        				edge from parent 
           					node[] {}
            				node[]  {}
    				}   
    				child {
        				node[probranch, label=above: $3/2$] {}    
        				edge from parent 
           					node[] {}
            				node[]  {}
    				}          
        edge from parent         
            node[] {}
            node[]  {}
    };
\end{tikzpicture}\\
\end{tabular}
\tikzstyle{level 1}=[level distance=6mm, sibling distance=4mm]
\tikzstyle{level 2}=[level distance=6mm, sibling distance=3mm]
\begin{tabular}{l r}
\begin{tikzpicture}[grow=right, sloped]
\draw  [-, rounded corners, thick] (.7,-.45) -- (.9,-.6) -- (.7,-.75);
\node[probranch] {}
    child {
        node[probranch] {}
        		child {
        				node[probranch] {}    
        				edge from parent 
           					node[] {}
            				node[]  {}
    				}    
        edge from parent 
            node[] {}
            node[]  {}
    }
    child {
        node[probranch] {}
        		child {
        				node[probranch] {}    
        				edge from parent 
           					node[] {}
            				node[]  {}
    				}    
        edge from parent 
            node[] {}
            node[]  {}
    }
    child {
        node[probranch] {}
        		child {
        				node[probranch] {}    
        				edge from parent 
           					node[] {}
            				node[]  {}
    				}    
        edge from parent 
            node[] {}
            node[]  {}
    }
    child {
        node[probranch] {}
        		child {
        				node[probranch] {}    
        				edge from parent 
           					node[] {}
            				node[]  {}
    				}    
        edge from parent 
            node[] {}
            node[]  {}
    }
    child {
        node[probranch, label=above: $1$] {}
        		child {
        				node[probranch] {}    
        				edge from parent 
           					node[] {}
            				node[]  {}
    				}   
    				child {
        				node[probranch, label=above: $3/2$] {}    
        				edge from parent 
           					node[] {}
            				node[]  {}
    				}          
        edge from parent         
            node[] {}
            node[]  {}
    };
\end{tikzpicture} \\
\end{tabular}
\tikzstyle{level 1}=[level distance=6mm, sibling distance=4mm]
\tikzstyle{level 2}=[level distance=6mm, sibling distance=3mm]
\begin{tabular}{l r}
\begin{tikzpicture}[grow=right, sloped]
\draw  [-, rounded corners, thick] (.7,-.45) -- (.9,-.6) -- (.7,-.75);
\draw[-,rounded corners,thick] (.7,.05)--(.9, .2)--(.7,.35);
\node[probranch] {}
    child {
        node[probranch] {}
        		child {
        				node[probranch] {}    
        				edge from parent 
           					node[] {}
            				node[]  {}
    				}    
        edge from parent 
            node[] {}
            node[]  {}
    }
    child {
        node[probranch] {}
        		child {
        				node[probranch] {}    
        				edge from parent 
           					node[] {}
            				node[]  {}
    				}    
        edge from parent 
            node[] {}
            node[]  {}
    }
    child {
        node[probranch] {}
        		child {
        				node[probranch] {}    
        				edge from parent 
           					node[] {}
            				node[]  {}
    				}    
        edge from parent 
            node[] {}
            node[]  {}
    }
    child {
        node[probranch] {}
        		child {
        				node[probranch] {}    
        				edge from parent 
           					node[] {}
            				node[]  {}
    				}    
        edge from parent 
            node[] {}
            node[]  {}
    }
    child {
        node[probranch, label=above: $1$] {}
        		child {
        				node[probranch] {}    
        				edge from parent 
           					node[] {}
            				node[]  {}
    				}   
    				child {
        				node[probranch, label=above: $3/2$] {}    
        				edge from parent 
           					node[] {}
            				node[]  {}
    				}          
        edge from parent         
            node[] {}
            node[]  {}
    };
\end{tikzpicture} \\
\end{tabular}
\tikzstyle{level 1}=[level distance=9mm, sibling distance=6mm]
\tikzstyle{level 2}=[level distance=6mm, sibling distance=3mm]
\begin{tabular}{l r}
\begin{tikzpicture}[grow=right, sloped]
\node[probranch] {}
    child {
        node[probranch] {}
        		child {
        				node[probranch] {}    
        				edge from parent 
           					node[] {}
            				node[]  {}
    				}    
    				child {
        				node[probranch] {}    
        				edge from parent 
           					node[] {}
            				node[]  {}
    				}      
        edge from parent 
            node[] {}
            node[]  {}
    }
    child {
        node[probranch] {}
        		child {
        				node[probranch] {}    
        				edge from parent 
           					node[] {}
            				node[]  {}
    				}
    				child {
        				node[probranch] {}    
        				edge from parent 
           					node[] {}
            				node[]  {}
    				}     
        edge from parent 
            node[] {}
            node[]  {}
    }
    child {
        node[probranch, label=above: $1$] {}
        		child {
        				node[probranch] {}    
        				edge from parent 
           					node[] {}
            				node[]  {}
    				} 
    				child {
        				node[probranch, label=above: $3/2$] {}    
        				edge from parent 
           					node[] {}
            				node[]  {}
    				}          
        edge from parent         
            node[] {}
            node[]  {}
    };
\end{tikzpicture}\\
\end{tabular}
\tikzstyle{level 1}=[level distance=9mm, sibling distance=6mm]
\tikzstyle{level 2}=[level distance=6mm, sibling distance=3mm]
\begin{tabular}{l r}
\begin{tikzpicture}[grow=right, sloped]
\draw  [-, rounded corners, thick] (1,-.1) -- (1.2,-.3) -- (1,-.5);
\node[probranch] {}
    child {
        node[probranch] {}
        		child {
        				node[probranch] {}    
        				edge from parent 
           					node[] {}
            				node[]  {}
    				}    
    				child {
        				node[probranch] {}    
        				edge from parent 
           					node[] {}
            				node[]  {}
    				}      
        edge from parent 
            node[] {}
            node[]  {}
    }
    child {
        node[probranch] {}
        		child {
        				node[probranch] {}    
        				edge from parent 
           					node[] {}
            				node[]  {}
    				}
    				child {
        				node[probranch] {}    
        				edge from parent 
           					node[] {}
            				node[]  {}
    				}     
        edge from parent 
            node[] {}
            node[]  {}
    }
    child {
        node[probranch, label=above: $1$] {}
        		child {
        				node[probranch] {}    
        				edge from parent 
           					node[] {}
            				node[]  {}
    				} 
    				child {
        				node[probranch, label=above: $3/2$] {}    
        				edge from parent 
           					node[] {}
            				node[]  {}
    				}          
        edge from parent         
            node[] {}
            node[]  {}
    };
\end{tikzpicture}\\
\end{tabular}
\tikzstyle{level 1}=[level distance=9mm, sibling distance=6mm]
\tikzstyle{level 2}=[level distance=6mm, sibling distance=3mm]
\begin{tabular}{l r}
\begin{tikzpicture}[grow=right, sloped]
\node[probranch] {}
    child {
        node[probranch] {}
        		child {
        				node[probranch] {}    
        				edge from parent 
           					node[] {}
            				node[]  {}
    				}    
        edge from parent 
            node[] {}
            node[]  {}
    }
    child {
        node[probranch] {}
        		child {
        				node[probranch] {}    
        				edge from parent 
           					node[] {}
            				node[]  {}
    				}    
        edge from parent 
            node[] {}
            node[]  {}
    }
    child {
        node[probranch] {}
        		child {
        				node[probranch] {}    
        				edge from parent 
           					node[] {}
            				node[]  {}
    				}
    				child {
        				node[probranch] {}    
        				edge from parent 
           					node[] {}
            				node[]  {}
    				}     
        edge from parent 
            node[] {}
            node[]  {}
    }
    child {
        node[probranch, label=above: $1$] {}
        		child {
        				node[probranch] {}    
        				edge from parent 
           					node[] {}
            				node[]  {}
    				}   
    				child {
        				node[probranch, label=above: $3/2$] {}    
        				edge from parent 
           					node[] {}
            				node[]  {}
    				}          
        edge from parent         
            node[] {}
            node[]  {}
    };
\end{tikzpicture}\\
\end{tabular}
\tikzstyle{level 1}=[level distance=9mm, sibling distance=6mm]
\tikzstyle{level 2}=[level distance=6mm, sibling distance=3mm]
\begin{tabular}{l r}
\begin{tikzpicture}[grow=right, sloped]
\draw  [-, rounded corners, thick] (1,-.35) -- (1.2,-.6) -- (1,-.85);
\node[probranch] {}
    child {
        node[probranch] {}
        		child {
        				node[probranch] {}    
        				edge from parent 
           					node[] {}
            				node[]  {}
    				}    
        edge from parent 
            node[] {}
            node[]  {}
    }
    child {
        node[probranch] {}
        		child {
        				node[probranch] {}    
        				edge from parent 
           					node[] {}
            				node[]  {}
    				}    
        edge from parent 
            node[] {}
            node[]  {}
    }
    child {
        node[probranch] {}
        		child {
        				node[probranch] {}    
        				edge from parent 
           					node[] {}
            				node[]  {}
    				}
    				child {
        				node[probranch] {}    
        				edge from parent 
           					node[] {}
            				node[]  {}
    				}     
        edge from parent 
            node[] {}
            node[]  {}
    }
    child {
        node[probranch, label=above: $1$] {}
        		child {
        				node[probranch] {}    
        				edge from parent 
           					node[] {}
            				node[]  {}
    				}   
    				child {
        				node[probranch, label=above: $3/2$] {}    
        				edge from parent 
           					node[] {}
            				node[]  {}
    				}          
        edge from parent         
            node[] {}
            node[]  {}
    };
\end{tikzpicture} \\
\end{tabular}
\end{center}
\noindent
Tangent Cones: \\ $y(y-x)(y+x)(y-ax)(y-bx)(y-cx) $ where $a, b, c\neq0,1,-1$ and $a\neq b \neq c$, \\ $y(y^2+x^2)(y-ax)(y-bx)(y-cx)$ where $a, b, c\neq0$ and $a\neq b \neq c$, \\ $y(y^2+x^2)(y^2+ax^2)(y-bx)$ where $a>0$, $a\neq1$, and $b\neq0$, \\ $y(y^2+x^2)^2(y-ax)$ where $a\neq0$ \\
\\
Newton Polygon 6:
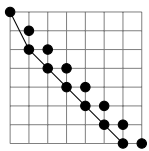
\begin{center}
\begin{tikzpicture}[scale=.25]
\draw[help lines] (0,0) grid (7,7);
\draw (0,7) -- (1,5);
\draw(1,5) -- (6,0);
\draw[fill] (0,7) circle [radius= 0.25];
\draw[fill] (1,6) circle [radius= 0.25];
\draw[fill] (2,5) circle [radius= 0.25];
\draw[fill] (3,4) circle [radius= 0.25];
\draw[fill] (4,3) circle [radius= 0.25];
\draw[fill] (5,2) circle [radius= 0.25];
\draw[fill] (6,1) circle [radius= 0.25];
\draw[fill] (7,0) circle [radius= 0.25];
\draw[fill] (6,0) circle [radius= 0.25];
\draw[fill](5,1) circle [radius = 0.25];
\draw[fill](4,2) circle [radius= 0.25];
\draw[fill](3,3) circle [radius = 0.25];
\draw[fill](2,4) circle[radius = 0.25];
\draw[fill](1,5) circle [radius = 0.25];
\end{tikzpicture}
\end{center} 
Families: \\ $y(y-x)(y+x)(y-ax)(y-bx)(y-cx)+dx^7+ex^6y+fx^5y^2+gx^4y^3+hx^3y^4+jx^2y^5+kxy^6+ly^7=0$ where $a, b, c\neq0,1,-1$ and $a\neq b \neq c$, \\ \\ $y(y^2+x^2)(y-ax)(y-bx)(y-cx)+dx^7+ex^6y+fx^5y^2+gx^4y^3+hx^3y^4+jx^2y^5+kxy^6+ly^7=0$  where $a, b, c\neq0$ and $a\neq b \neq c$, \\ \\ $y(y^2+x^2)(y^2+ax^2)(y-bx)+cx^7+dx^6y+ex^5y^2+fx^4y^3+gx^3y^4+hx^2y^5+jxy^6+ky^7=0$ where $a>0$, $a\neq1$, and $b\neq0$, \\ \\ $y(y^2+x^2)^2(y-ax)+bx^7+cx^6y+dx^5y^2+ex^4y^3+fx^3y^4+gx^2y^5+hxy^6+jy^7=0$ where $a\neq0$ \\
\\
Diagrams: 
\begin{center}
\tikzstyle{level 1}=[level distance=9mm, sibling distance=3mm]
\tikzstyle{level 2}=[level distance=6mm, sibling distance=6mm]
\begin{tabular}{l r}
\begin{tikzpicture}[grow=right, sloped]
\node[probranch] {}
    child {
        node[probranch] {}        
        edge from parent 
            node[] {}
            node[]  {}
    }
    child {
        node[probranch] {}        
        edge from parent 
            node[] {}
            node[]  {}
    }
    child {
        node[probranch] {}        
        edge from parent 
            node[] {}
            node[]  {}
    }
    child {
        node[probranch] {}        
        edge from parent 
            node[] {}
            node[]  {}
    }
    child {
        node[probranch] {}        
        edge from parent 
            node[] {}
            node[]  {}
    }
    child {
        node[probranch, label=above: $1$] {}        
        edge from parent         
            node[] {}
            node[]  {}
    };
\end{tikzpicture}\\
\end{tabular}
\tikzstyle{level 1}=[level distance=9mm, sibling distance=3mm]
\tikzstyle{level 2}=[level distance=6mm, sibling distance=6mm]
\begin{tabular}{l r}
\begin{tikzpicture}[grow=right, sloped]
\draw  [-, rounded corners, thick] (1,-.45) -- (1.2,-.6) -- (1,-.75);
\node[probranch] {}
    child {
        node[probranch] {}        
        edge from parent 
            node[] {}
            node[]  {}
    }
    child {
        node[probranch] {}        
        edge from parent 
            node[] {}
            node[]  {}
    }
    child {
        node[probranch] {}        
        edge from parent 
            node[] {}
            node[]  {}
    }
    child {
        node[probranch] {}        
        edge from parent 
            node[] {}
            node[]  {}
    }
    child {
        node[probranch] {}        
        edge from parent 
            node[] {}
            node[]  {}
    }
    child {
        node[probranch, label=above: $1$] {}        
        edge from parent         
            node[] {}
            node[]  {}
    };
\end{tikzpicture} \\
\end{tabular}
\tikzstyle{level 1}=[level distance=9mm, sibling distance=3mm]
\tikzstyle{level 2}=[level distance=6mm, sibling distance=6mm]
\begin{tabular}{l r}
\begin{tikzpicture}[grow=right, sloped]
\draw  [-, rounded corners, thick] (1,-.45) -- (1.2,-.6) -- (1,-.75);
\draw [-, rounded corners, thick] (1, -.15) -- (1.2, 0) -- (1, .15);
\node[probranch] {}
    child {
        node[probranch] {}        
        edge from parent 
            node[] {}
            node[]  {}
    }
    child {
        node[probranch] {}        
        edge from parent 
            node[] {}
            node[]  {}
    }
    child {
        node[probranch] {}        
        edge from parent 
            node[] {}
            node[]  {}
    }
    child {
        node[probranch] {}        
        edge from parent 
            node[] {}
            node[]  {}
    }
    child {
        node[probranch] {}        
        edge from parent 
            node[] {}
            node[]  {}
    }
    child {
        node[probranch, label=above: $1$] {}        
        edge from parent         
            node[] {}
            node[]  {}
    };
\end{tikzpicture} \\
\end{tabular}
\tikzstyle{level 1}=[level distance=9mm, sibling distance=6mm]
\tikzstyle{level 2}=[level distance=6mm, sibling distance=3mm]
\begin{tabular}{l r}
\begin{tikzpicture}[grow=right, sloped]
\draw  [-, rounded corners, thick] (1,-.35) -- (1.2,-.6) -- (1,-.85);
\node[probranch] {}
    child {
        node[probranch] {}
        		child {
        				node[probranch] {}    
        				edge from parent 
           					node[] {}
            				node[]  {}
    				}    
        edge from parent 
            node[] {}
            node[]  {}
    }
    child {
        node[probranch] {}
        		child {
        				node[probranch] {}    
        				edge from parent 
           					node[] {}
            				node[]  {}
    				}    
        edge from parent 
            node[] {}
            node[]  {}
    }
    child {
        node[probranch] {}
        		child {
        				node[probranch] {}    
        				edge from parent 
           					node[] {}
            				node[]  {}
    				}
    				child {
        				node[probranch] {}    
        				edge from parent 
           					node[] {}
            				node[]  {}
    				}     
        edge from parent 
            node[] {}
            node[]  {}
    }
    child {
        node[probranch, label=above: $1$] {}
        		child {
        				node[probranch] {}    
        				edge from parent 
           					node[] {}
            				node[]  {}
    				}   
    				child {
        				node[probranch, label=above: $3/2$] {}    
        				edge from parent 
           					node[] {}
            				node[]  {}
    				}          
        edge from parent         
            node[] {}
            node[]  {}
    };
\end{tikzpicture} \\
\end{tabular}%
\end{center} 
Families: \\ $(y^2+x^2)(y^2+ax^2)(y^2+bx^2)+cx^7+dx^6y+ex^5y^2+fx^4y^3+gx^3y^4+hx^2y^5+jxy^6+ky^7=0$ where $a,b>0$, $a, b\neq1$, $a\neq b$, \\ \\ $(y^2+x^2)^2(y^2+ax^2)+cx^7+dx^6y+ex^5y^2+fx^4y^3+gx^3y^4+hx^2y^5+jxy^6+ky^7=0 $ where $a>0$ and $a\neq1$ \\ \\ $(y^2+x^2)^3+cx^7+dx^6y+ex^5y^2+fx^4y^3+gx^3y^4+hx^2y^5+jxy^6+ky^7=0$
\\
Diagrams: 
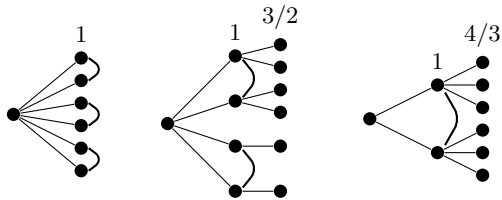
\begin{center}
\tikzstyle{level 1}=[level distance=9mm, sibling distance=3mm]
\tikzstyle{level 2}=[level distance=6mm, sibling distance=6mm]
\begin{tabular}{l r}
\begin{tikzpicture}[grow=right, sloped]
\draw  [-, rounded corners, thick] (1,-.45) -- (1.2,-.6) -- (1,-.75);
\draw [-, rounded corners, thick] (1, -.15) -- (1.2, 0) -- (1, .15);
\draw [-, rounded corners, thick] (1, .45) -- (1.2, .6) -- (1, .75);
\node[probranch] {}
    child {
        node[probranch] {}        
        edge from parent 
            node[] {}
            node[]  {}
    }
    child {
        node[probranch] {}        
        edge from parent 
            node[] {}
            node[]  {}
    }
    child {
        node[probranch] {}        
        edge from parent 
            node[] {}
            node[]  {}
    }
    child {
        node[probranch] {}        
        edge from parent 
            node[] {}
            node[]  {}
    }
    child {
        node[probranch] {}        
        edge from parent 
            node[] {}
            node[]  {}
    }
    child {
        node[probranch, label=above: $1$] {}        
        edge from parent         
            node[] {}
            node[]  {}
    };
\end{tikzpicture}\\
\end{tabular}
\tikzstyle{level 1}=[level distance=9mm, sibling distance=6mm]
\tikzstyle{level 2}=[level distance=6mm, sibling distance=3mm]
\begin{tabular}{l r}
\begin{tikzpicture}[grow=right, sloped]
\draw  [-, rounded corners, thick] (1,-.35) -- (1.2,-.6) -- (1,-.85);
\draw  [-, rounded corners, thick] (1,.35) -- (1.2,.6) -- (1,.85);
\node[probranch] {}
    child {
        node[probranch] {}
        		child {
        				node[probranch] {}    
        				edge from parent 
           					node[] {}
            				node[]  {}
    				}    
        edge from parent 
            node[] {}
            node[]  {}
    }
    child {
        node[probranch] {}
        		child {
        				node[probranch] {}    
        				edge from parent 
           					node[] {}
            				node[]  {}
    				}    
        edge from parent 
            node[] {}
            node[]  {}
    }
    child {
        node[probranch] {}
        		child {
        				node[probranch] {}    
        				edge from parent 
           					node[] {}
            				node[]  {}
    				}
    				child {
        				node[probranch] {}    
        				edge from parent 
           					node[] {}
            				node[]  {}
    				}     
        edge from parent 
            node[] {}
            node[]  {}
    }
    child {
        node[probranch, label=above: $1$] {}
        		child {
        				node[probranch] {}    
        				edge from parent 
           					node[] {}
            				node[]  {}
    				}   
    				child {
        				node[probranch, label=above: $3/2$] {}    
        				edge from parent 
           					node[] {}
            				node[]  {}
    				}          
        edge from parent         
            node[] {}
            node[]  {}
    };
\end{tikzpicture} \\
\end{tabular}
\tikzstyle{level 1}=[level distance=9mm, sibling distance=9mm]
\tikzstyle{level 2}=[level distance=6mm, sibling distance=3mm]
\begin{tabular}{l r}
\begin{tikzpicture}[grow=right, sloped]
\draw  [-, rounded corners, thick] (1,-.35) -- (1.2,0) -- (1,.35);
\node[probranch] {}
    child {
        node[probranch] {}
        		child {
        				node[probranch] {}    
        				edge from parent 
           					node[] {}
            				node[]  {}
    				}
    				child {
        				node[probranch] {}    
        				edge from parent 
           					node[] {}
            				node[]  {}
    				}
    				child {
        				node[probranch] {}    
        				edge from parent 
           					node[] {}
            				node[]  {}
    				}    
        edge from parent 
            node[] {}
            node[]  {}
    }
    child {
        node[probranch, label=above: $1$] {}
        		child {
        				node[probranch] {}    
        				edge from parent 
           					node[] {}
            				node[]  {}
    				}
    				child {
        				node[probranch] {}    
        				edge from parent 
           					node[] {}
            				node[]  {}
    				}      
    				child {
        				node[probranch, label=above: $4/3$] {}    
        				edge from parent 
           					node[] {}
            				node[]  {}
    				}          
        edge from parent         
            node[] {}
            node[]  {}
    };
\end{tikzpicture} \\
\end{tabular}
\end{center}

\section{Real Singular Points of Multiplicity 5}

\noindent Tangent Cone: $y^5$ \\ 
\\
Newton Polygons 8, 9: \\
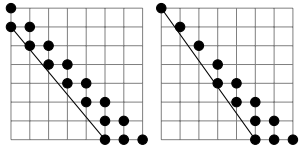
\begin{center}
%
\end{center} 
Families: \\ $y^5+ax^6+bx^5y+cx^4y^2+dx^3y^3+ex^2y^4+fxy^5+gy^6+hx^7+jx^6y+kx^5y^2+lx^4y^3+mx^3y^4+nx^2y^5+pxy^6+qy^7=0$, \\ \\ $y^5+ax^3y^3+bx^2y^4+cxy^5+dy^6+ex^7+fx^6y+gx^5y^2+hx^4y^3+jx^3y^4+kx^2y^5+lxy^6+my^7=0$ \\ \\
Diagrams: 
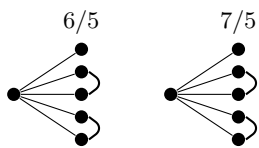
\begin{center}

\tikzstyle{level 1}=[level distance=9mm, sibling distance=3mm]
\tikzstyle{level 2}=[level distance=6mm, sibling distance=6mm]
\begin{tabular}{l r}
\begin{tikzpicture}[grow=right, sloped]
\draw  [-, rounded corners, thick] (1,-.3) -- (1.2,-.45) -- (1,-.6);
\draw [-, rounded corners, thick] (1, 0) -- (1.2, .15) -- (1, .3);
\node[probranch] {}
    child {
        node[probranch] {}        
        edge from parent 
            node[] {}
            node[]  {}
    }
    child {
        node[probranch] {}        
        edge from parent 
            node[] {}
            node[]  {}
    }
    child {
        node[probranch] {}        
        edge from parent 
            node[] {}
            node[]  {}
    }
    child {
        node[probranch] {}        
        edge from parent 
            node[] {}
            node[]  {}
    }
    child {
        node[probranch, label=above: $6/5$] {}        
        edge from parent         
            node[] {}
            node[]  {}
    };
\end{tikzpicture}\\
\end{tabular}
\tikzstyle{level 1}=[level distance=9mm, sibling distance=3mm]
\tikzstyle{level 2}=[level distance=6mm, sibling distance=6mm]
\begin{tabular}{l r}
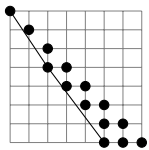
\begin{tikzpicture}[grow=right, sloped]
\draw  [-, rounded corners, thick] (1,-.3) -- (1.2,-.45) -- (1,-.6);
\draw [-, rounded corners, thick] (1, 0) -- (1.2, .15) -- (1, .3);
\node[probranch] {}
    child {
        node[probranch] {}        
        edge from parent 
            node[] {}
            node[]  {}
    }
    child {
        node[probranch] {}        
        edge from parent 
            node[] {}
            node[]  {}
    }
    child {
        node[probranch] {}        
        edge from parent 
            node[] {}
            node[]  {}
    }
    child {
        node[probranch] {}        
        edge from parent 
            node[] {}
            node[]  {}
    }
    child {
        node[probranch, label=above: $7/5$] {}        
        edge from parent         
            node[] {}
            node[]  {}
    };
\end{tikzpicture}\\
\end{tabular}%
\end{center}
Family: \\ $y^5+ax^4y^2+bx^3y^3+cx^2y^4+dxy^5+ey^6+fx^7+gx^6y+hx^5y^2+jx^4y^3+kx^3y^4+lx^2y^5+mxy^6+ny^7=0$ \\ \\
Diagram: 
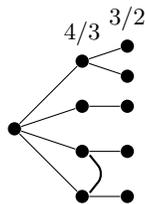
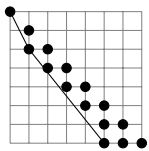
\begin{center}
\tikzstyle{level 1}=[level distance=9mm, sibling distance=6mm]
\tikzstyle{level 2}=[level distance=6mm, sibling distance=4mm]
\begin{tabular}{l r}
\begin{tikzpicture}[grow=right, sloped]
\draw  [-, rounded corners, thick] (1,-.35) -- (1.2,-.6) -- (1,-.85);
\node[probranch] {}
    child {
        node[probranch] {} 
           child {
        				node[probranch] {}    
        				edge from parent 
           					node[] {}
            				node[]  {}
            }     
        edge from parent 
            node[] {}
            node[]  {}
    }
    child {
        node[probranch] {}
        		child {
        				node[probranch] {}    
        				edge from parent 
           					node[] {}
            				node[]  {}
    				}    
        edge from parent 
            node[] {}
            node[]  {}
    }
    child {
        node[probranch] {}
        		child {
        				node[probranch] {}    
        				edge from parent 
           					node[] {}
            				node[]  {}
    				}    
        edge from parent 
            node[] {}
            node[]  {}
    }
    child {
        node[probranch, label=above: $4/3$] {}
        		child {
        				node[probranch] {}    
        				edge from parent 
           					node[] {}
            				node[]  {}
    				}
    				child {
        				node[probranch, label=above: $3/2$] {}    
        				edge from parent 
           					node[] {}
            				node[]  {}
    				}          
        edge from parent         
            node[] {}
            node[]  {}
    };
\end{tikzpicture}\\
\end{tabular}%
\end{center}
Family: \\ $y^5+ax^5y+bx^4y^2+cx^3y^3+dx^2y^4+exy^5+fy^6+gx^7+hx^6y+jx^5y^2+kx^4y^3+lx^3y^4+mx^2y^5+nxy^6+py^7=0$ \\ \\
Diagram: \begin{center}
\tikzstyle{level 1}=[level distance=9mm, sibling distance=3mm]
\tikzstyle{level 2}=[level distance=6mm, sibling distance=6mm]
\begin{tabular}{l r}
\begin{tikzpicture}[grow=right, sloped]
\draw  [-, rounded corners, thick] (1,-.3) -- (1.2,-.45) -- (1,-.6);
\node[probranch] {}
    child {
        node[probranch] {}        
        edge from parent 
            node[] {}
            node[]  {}
    }
    child {
        node[probranch] {}        
        edge from parent 
            node[] {}
            node[]  {}
    }
    child {
        node[probranch] {}        
        edge from parent 
            node[] {}
            node[]  {}
    }
    child {
        node[probranch] {}        
        edge from parent 
            node[] {}
            node[]  {}
    }
    child {
        node[probranch, label=above: $5/4$] {}        
        edge from parent         
            node[] {}
            node[]  {}
    };
\end{tikzpicture}\\
\end{tabular}%
\end{center}
Family: \\ $y^4(y-x)+ax^4y^2+bx^3y^3+cx^2y^4+dxy^5+ey^6+fx^7+gx^6y+hx^5y^2+jx^4y^3+kx^3y^4+lx^2y^5+mxy^6+ny^7=0$ \\ \\
Diagrams: 
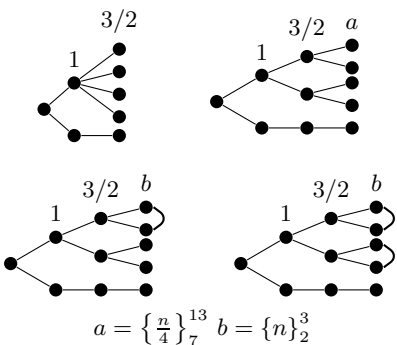
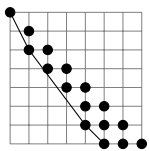
\begin{center}
\tikzstyle{level 1}=[level distance=4mm, sibling distance=7mm]
\tikzstyle{level 2}=[level distance=6mm, sibling distance=3mm]
\begin{tabular}{l r}
\begin{tikzpicture}[grow=right, sloped]
\node[probranch] {}
    child {
        node[probranch] {}
        		child {
        				node[probranch] {}    
        				edge from parent 
           					node[] {}
            				node[]  {}
    				}    
        edge from parent 
            node[] {}
            node[]  {}
    }
    child {
        node[probranch, label=above: $1$] {}
            child {
        				node[probranch] {}    
        				edge from parent 
           					node[] {}
            				node[]  {}
    				}
        		child {
        				node[probranch] {}    
        				edge from parent 
           					node[] {}
            				node[]  {}
    				}
    				child {
        				node[probranch] {}    
        				edge from parent 
           					node[] {}
            				node[]  {}
    				}    
    				child {
        				node[probranch, label=above: $3/2$] {}    
        				edge from parent 
           					node[] {}
            				node[]  {}
    				}          
        edge from parent         
            node[] {}
            node[]  {}
    };
\end{tikzpicture}\\
\end{tabular}
\tikzstyle{level 1}=[level distance=6mm, sibling distance=7mm]
\tikzstyle{level 2}=[level distance=6mm, sibling distance=5mm]
\tikzstyle{level 3}=[level distance=6mm, sibling distance=3mm]
\begin{tabular}{l r}
\begin{tikzpicture}[grow=right, sloped]
\node[probranch] {}
    child {
        node[probranch] {}
        		child {
        				node[probranch] {}
        						child {
        								node[probranch] {}    
        						edge from parent 
           							node[] {}
            						node[]  {}
    								}   
        				edge from parent 
           					node[] {}
            				node[]  {}
    				}   
        edge from parent 
            node[] {}
            node[]  {}
    }
    child {
        node[probranch, label=above: $1$] {}
            child {
        				node[probranch] {}
        						child {
        								node[probranch] {}    
        						edge from parent 
           							node[] {}
            						node[]  {}
    								}
    								child {
        								node[probranch] {}    
        						edge from parent 
           							node[] {}
            						node[]  {}
    								}     
        				edge from parent 
           					node[] {}
            				node[]  {}
    				} 
    				child {
        				node[probranch, label=above: $3/2$] {}
        						child {
        								node[probranch] {}    
        						edge from parent 
           							node[] {}
            						node[]  {}
    								} 
        						child {
        								node[probranch, label=above: $a$] {}    
        						edge from parent 
           							node[] {}
            						node[]  {}
    								}    
        				edge from parent 
           					node[] {}
            				node[]  {}
    				}          
        edge from parent         
            node[] {}
            node[]  {}
    };
\end{tikzpicture}\\
\end{tabular}
\end{center}
\begin{center}
\tikzstyle{level 1}=[level distance=6mm, sibling distance=7mm]
\tikzstyle{level 2}=[level distance=6mm, sibling distance=5mm]
\tikzstyle{level 3}=[level distance=6mm, sibling distance=3mm]
\begin{tabular}{l r}
\begin{tikzpicture}[grow=right, sloped]
\draw  [-, rounded corners, thick] (1.9,.45) -- (2.1,.6) -- (1.9,.75);
\node[probranch] {}
    child {
        node[probranch] {}
        		child {
        				node[probranch] {}
        						child {
        								node[probranch] {}    
        						edge from parent 
           							node[] {}
            						node[]  {}
    								}   
        				edge from parent 
           					node[] {}
            				node[]  {}
    				}   
        edge from parent 
            node[] {}
            node[]  {}
    }
    child {
        node[probranch, label=above: $1$] {}
            child {
        				node[probranch] {}
        						child {
        								node[probranch] {}    
        						edge from parent 
           							node[] {}
            						node[]  {}
    								}
    								child {
        								node[probranch] {}    
        						edge from parent 
           							node[] {}
            						node[]  {}
    								}     
        				edge from parent 
           					node[] {}
            				node[]  {}
    				} 
    				child {
        				node[probranch, label=above: $3/2$] {}
        						child {
        								node[probranch] {}    
        						edge from parent 
           							node[] {}
            						node[]  {}
    								} 
        						child {
        								node[probranch, label=above: $b$] {}    
        						edge from parent 
           							node[] {}
            						node[]  {}
    								}    
        				edge from parent 
           					node[] {}
            				node[]  {}
    				}          
        edge from parent         
            node[] {}
            node[]  {}
    };
\end{tikzpicture}\\
\end{tabular}
\tikzstyle{level 1}=[level distance=6mm, sibling distance=7mm]
\tikzstyle{level 2}=[level distance=6mm, sibling distance=5mm]
\tikzstyle{level 3}=[level distance=6mm, sibling distance=3mm]
\begin{tabular}{l r}
\begin{tikzpicture}[grow=right, sloped]
\draw  [-, rounded corners, thick] (1.9,.45) -- (2.1,.6) -- (1.9,.75);
\draw [-, rounded corners, thick] (1.9, .25) -- (2.1, .1) -- (1.9, -.05);
\node[probranch] {}
    child {
        node[probranch] {}
        		child {
        				node[probranch] {}
        						child {
        								node[probranch] {}    
        						edge from parent 
           							node[] {}
            						node[]  {}
    								}   
        				edge from parent 
           					node[] {}
            				node[]  {}
    				}   
        edge from parent 
            node[] {}
            node[]  {}
    }
    child {
        node[probranch, label=above: $1$] {}
            child {
        				node[probranch] {}
        						child {
        								node[probranch] {}    
        						edge from parent 
           							node[] {}
            						node[]  {}
    								}
    								child {
        								node[probranch] {}    
        						edge from parent 
           							node[] {}
            						node[]  {}
    								}     
        				edge from parent 
           					node[] {}
            				node[]  {}
    				} 
    				child {
        				node[probranch, label=above: $3/2$] {}
        						child {
        								node[probranch] {}    
        						edge from parent 
           							node[] {}
            						node[]  {}
    								} 
        						child {
        								node[probranch, label=above: $b$] {}    
        						edge from parent 
           							node[] {}
            						node[]  {}
    								}    
        				edge from parent 
           					node[] {}
            				node[]  {}
    				}          
        edge from parent         
            node[] {}
            node[]  {}
    };
\end{tikzpicture}\\
\end{tabular}%
\end{center}
Family: \\ $y^4(y-x)+ax^5y+bx^4y^2+cx^3y^3+dx^2y^4+exy^5+fy^6+gx^7+hx^6y+jx^5y^2+kx^4y^3+lx^3y^4+mx^2y^5+nxy^6+py^7=0$ \\ \\
Diagram: 
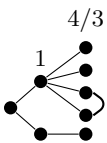
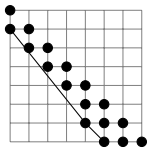
\begin{center}
\tikzstyle{level 1}=[level distance=4mm, sibling distance=7mm]
\tikzstyle{level 2}=[level distance=6mm, sibling distance=3mm]
\begin{tabular}{l r}
\begin{tikzpicture}[grow=right, sloped]
\draw  [-, rounded corners, thick] (1.1,-.1) -- (1.3,.05) -- (1.1,.2);
\node[probranch] {}
    child {
        node[probranch] {}
        		child {
        				node[probranch] {}    
        				edge from parent 
           					node[] {}
            				node[]  {}
    				}    
        edge from parent 
            node[] {}
            node[]  {}
    }
    child {
        node[probranch, label=above: $1$] {}
            child {
        				node[probranch] {}    
        				edge from parent 
           					node[] {}
            				node[]  {}
    				}
        		child {
        				node[probranch] {}    
        				edge from parent 
           					node[] {}
            				node[]  {}
    				}
    				child {
        				node[probranch] {}    
        				edge from parent 
           					node[] {}
            				node[]  {}
    				}    
    				child {
        				node[probranch, label=above: $4/3$] {}    
        				edge from parent 
           					node[] {}
            				node[]  {}
    				}          
        edge from parent         
            node[] {}
            node[]  {}
    };
\end{tikzpicture}\\
\end{tabular}%
\end{center}
Family: $y^4(y-x)+ax^6+bx^5y+cx^4y^2+dx^3y^3+ex^2y^4+fxy^5+gy^6+hx^7+jx^6y+kx^5y^2+lx^4y^3+mx^3y^4+nx^2y^5+pxy^6+qy^7=0$ \\ \\
Diagram: 
\begin{center}
\tikzstyle{level 1}=[level distance=4mm, sibling distance=7mm]
\tikzstyle{level 2}=[level distance=6mm, sibling distance=3mm]
\begin{tabular}{l r}
\begin{tikzpicture}[grow=right, sloped]
\draw  [-, rounded corners, thick] (1.1,-.1) -- (1.3,.05) -- (1.1,.2);
\node[probranch] {}
    child {
        node[probranch] {}
        		child {
        				node[probranch] {}    
        				edge from parent 
           					node[] {}
            				node[]  {}
    				}    
        edge from parent 
            node[] {}
            node[]  {}
    }
    child {
        node[probranch, label=above: $1$] {}
            child {
        				node[probranch] {}    
        				edge from parent 
           					node[] {}
            				node[]  {}
    				}
        		child {
        				node[probranch] {}    
        				edge from parent 
           					node[] {}
            				node[]  {}
    				}
    				child {
        				node[probranch] {}    
        				edge from parent 
           					node[] {}
            				node[]  {}
    				}    
    				child {
        				node[probranch, label=above: $5/4$] {}    
        				edge from parent 
           					node[] {}
            				node[]  {}
    				}          
        edge from parent         
            node[] {}
            node[]  {}
    };
\end{tikzpicture}\\
\end{tabular}%
\end{center}
Families: \\ $y^3(y-x)^2+ax^4y^2+bx^3y^3+cx^2y^4+dxy^5+ey^6+fx^7+gx^6y+hx^5y^2+jx^4y^3+kx^3y^4+lx^2y^5+mxy^6+ny^7=0$, \\ \\ $y^3(y-x)(y-2x)+ax^4y^2+bx^3y^3+cx^2y^4+dxy^5+ey^6+fx^7+gx^6y+hx^5y^2+jx^4y^3+kx^3y^4+lx^2y^5+mxy^6+ny^7=0$, \\ \\ $y^3(y^2+x^2)+ax^4y^2+bx^3y^3+cx^2y^4+dxy^5+ey^6+fx^7+gx^6y+hx^5y^2+jx^4y^3+kx^3y^4+lx^2y^5+mxy^6+ny^7=0$, \\ \\ $y^3(y-x)^2+ax^6+bx^5y+cx^4y^2+dx^3y^3+ex^2y^4+fxy^5+gy^6+hx^7+jx^6y+kx^5y^2+lx^4y^3+mx^3y^4+nx^2y^5+pxy^6+qy^7=0$, \\ \\ $y^3(y-x)(y-2x)+ax^6+bx^5y+cx^4y^2+dx^3y^3+ex^2y^4+fxy^5+gy^6+hx^7+jx^6y+kx^5y^2+lx^4y^3+mx^3y^4+nx^2y^5+pxy^6+qy^7=0$, \\ \\ $y^3(y^2+x^2)+ax^6+bx^5y+cx^4y^2+dx^3y^3+ex^2y^4+fxy^5+gy^6+hx^7+jx^6y+kx^5y^2+lx^4y^3+mx^3y^4+nx^2y^5+pxy^6+qy^7=0$ \\ \\
Diagrams: 
\begin{center}
\tikzstyle{level 1}=[level distance=9mm, sibling distance=9mm]
\tikzstyle{level 2}=[level distance=6mm, sibling distance=3mm]
\begin{tabular}{l r}
\begin{tikzpicture}[grow=right, sloped]
\draw  [-, rounded corners, thick] (1.6,1.25) -- (1.8,1.05) -- (1.6,.9);
\node[probranch] {}
    child {
        node[probranch] {}
        		child {
        				node[probranch] {}    
        				edge from parent 
           					node[] {}
            				node[]  {}
    				}    
        edge from parent 
            node[] {}
            node[]  {}
    }
    child {
        node[probranch] {}
        		child {
        				node[probranch] {}    
        				edge from parent 
           					node[] {}
            				node[]  {}
    				}    
        edge from parent 
            node[] {}
            node[]  {}
    }
    child {
        node[probranch, label=above: $1$] {}
        		child {
        				node[probranch] {}    
        				edge from parent 
           					node[] {}
            				node[]  {}
    				}
    				child {
        				node[probranch] {}    
        				edge from parent 
           					node[] {}
            				node[]  {}
    				}    
    				child {
        				node[probranch, label=above: $5/3$] {}    
        				edge from parent 
           					node[] {}
            				node[]  {}
    				}          
        edge from parent         
            node[] {}
            node[]  {}
    };
\end{tikzpicture} \\
\end{tabular}
\tikzstyle{level 1}=[level distance=9mm, sibling distance=9mm]
\tikzstyle{level 2}=[level distance=6mm, sibling distance=3mm]
\begin{tabular}{l r}
\begin{tikzpicture}[grow=right, sloped]
\draw  [-, rounded corners, thick] (1.6,1.25) -- (1.8,1.05) -- (1.6,.9);
\draw [-, rounded corners, thick] (1, -.05) -- (1.2, -.45) -- (1, -.85);
\node[probranch] {}
    child {
        node[probranch] {}
        		child {
        				node[probranch] {}    
        				edge from parent 
           					node[] {}
            				node[]  {}
    				}    
        edge from parent 
            node[] {}
            node[]  {}
    }
    child {
        node[probranch] {}
        		child {
        				node[probranch] {}    
        				edge from parent 
           					node[] {}
            				node[]  {}
    				}    
        edge from parent 
            node[] {}
            node[]  {}
    }
    child {
        node[probranch, label=above: $1$] {}
        		child {
        				node[probranch] {}    
        				edge from parent 
           					node[] {}
            				node[]  {}
    				}
    				child {
        				node[probranch] {}    
        				edge from parent 
           					node[] {}
            				node[]  {}
    				}    
    				child {
        				node[probranch, label=above: $5/3$] {}    
        				edge from parent 
           					node[] {}
            				node[]  {}
    				}          
        edge from parent         
            node[] {}
            node[]  {}
    };
\end{tikzpicture} \\
\end{tabular}
\tikzstyle{level 1}=[level distance=6mm, sibling distance=7mm]
\tikzstyle{level 2}=[level distance=6mm, sibling distance=4mm]
\tikzstyle{level 3}=[level distance=6mm, sibling distance=2mm]
\begin{tabular}{l r}
\begin{tikzpicture}[grow=right, sloped]
\draw  [-, rounded corners, thick] (1.9,-.35) -- (2.1,-.45) -- (1.9,-.55);
\node[probranch] {}
    child {
        node[probranch] {}
        		child {
        				node[probranch] {}
        						child {
        								node[probranch] {}    
        						edge from parent 
           							node[] {}
            						node[]  {}
    								}
    								child {
        								node[probranch] {}    
        						edge from parent 
           							node[] {}
            						node[]  {}
    								}      
							child {
        								node[probranch] {}    
        						edge from parent 
           							node[] {}
            						node[]  {}
    								}
        				edge from parent 
           					node[] {}
            				node[]  {}
    				}   
        edge from parent 
            node[] {}
            node[]  {}
    }
    child {
        node[probranch, label=above: $1$] {}
            child {
        				node[probranch] {}
        						child {
        								node[probranch] {}    
        						edge from parent 
           							node[] {}
            						node[]  {}
    								}    
        				edge from parent 
           					node[] {}
            				node[]  {}
    				} 
    				child {
        				node[probranch, label=above: $3/2$] {}
        						child {
        								node[probranch, label=above: $5/3$] {}    
        						edge from parent 
           							node[] {}
            						node[]  {}
    								}    
        				edge from parent 
           					node[] {}
            				node[]  {}
    				}          
        edge from parent         
            node[] {}
            node[]  {}
    };
\end{tikzpicture} \\
\end{tabular}
\tikzstyle{level 1}=[level distance=6mm, sibling distance=7mm]
\tikzstyle{level 2}=[level distance=6mm, sibling distance=4mm]
\tikzstyle{level 3}=[level distance=6mm, sibling distance=2mm]
\begin{tabular}{l r}
\begin{tikzpicture}[grow=right, sloped]
\draw  [-, rounded corners, thick] (1.3,0) -- (1.5,.15) -- (1.3,.3);
\node[probranch] {}
    child {
        node[probranch] {}
        		child {
        				node[probranch] {}
        						child {
        								node[probranch] {}    
        						edge from parent 
           							node[] {}
            						node[]  {}
    								}
    								child {
        								node[probranch] {}    
        						edge from parent 
           							node[] {}
            						node[]  {}
    								}      
        				edge from parent 
           					node[] {}
            				node[]  {}
    				}   
        edge from parent 
            node[] {}
            node[]  {}
    }
    child {
        node[probranch, label=above: $1$] {}
            child {
        				node[probranch] {}
        						child {
        								node[probranch] {}    
        						edge from parent 
           							node[] {}
            						node[]  {}
    								}    
        				edge from parent 
           					node[] {}
            				node[]  {}
    				} 
        		child {
        				node[probranch] {}
        						child {
        								node[probranch] {}    
        						edge from parent 
           							node[] {}
            						node[]  {}
    								}    
        				edge from parent 
           					node[] {}
            				node[]  {}
    				} 
    				child {
        				node[probranch, label=above: $5/3$] {}
        						child {
        								node[probranch, label=above: $a$] {}    
        						edge from parent 
           							node[] {}
            						node[]  {}
    								}    
        				edge from parent 
           					node[] {}
            				node[]  {}
    				}          
        edge from parent         
            node[] {}
            node[]  {}
    };
\end{tikzpicture}\\
\end{tabular}
\tikzstyle{level 1}=[level distance=6mm, sibling distance=7mm]
\tikzstyle{level 2}=[level distance=6mm, sibling distance=4mm]
\tikzstyle{level 3}=[level distance=6mm, sibling distance=2mm]
\begin{tabular}{l r}
\begin{tikzpicture}[grow=right, sloped]
\draw  [-, rounded corners, thick] (1.3,0) -- (1.5,.15) -- (1.3,.3);
\draw  [-, rounded corners, thick] (1.9,-.25) -- (2.1,-.35) -- (1.9,-.45);
\node[probranch] {}
    child {
        node[probranch] {}
        		child {
        				node[probranch] {}
        						child {
        								node[probranch] {}    
        						edge from parent 
           							node[] {}
            						node[]  {}
    								}
    								child {
        								node[probranch] {}    
        						edge from parent 
           							node[] {}
            						node[]  {}
    								}      
        				edge from parent 
           					node[] {}
            				node[]  {}
    				}   
        edge from parent 
            node[] {}
            node[]  {}
    }
    child {
        node[probranch, label=above: $1$] {}
            child {
        				node[probranch] {}
        						child {
        								node[probranch] {}    
        						edge from parent 
           							node[] {}
            						node[]  {}
    								}    
        				edge from parent 
           					node[] {}
            				node[]  {}
    				} 
        		child {
        				node[probranch] {}
        						child {
        								node[probranch] {}    
        						edge from parent 
           							node[] {}
            						node[]  {}
    								}    
        				edge from parent 
           					node[] {}
            				node[]  {}
    				} 
    				child {
        				node[probranch, label=above: $5/3$] {}
        						child {
        								node[probranch, label=above: $b$] {}    
        						edge from parent 
           							node[] {}
            						node[]  {}
    								}    
        				edge from parent 
           					node[] {}
            				node[]  {}
    				}          
        edge from parent         
            node[] {}
            node[]  {}
    };
\end{tikzpicture} \\
\end{tabular}
\tikzstyle{level 1}=[level distance=9mm, sibling distance=9mm]
\tikzstyle{level 2}=[level distance=6mm, sibling distance=3mm]
\begin{tabular}{l r}
\begin{tikzpicture}[grow=right, sloped]
\draw  [-, rounded corners, thick] (1.6,.9) -- (1.8,1.05) -- (1.6,1.2);
\node[probranch] {}
    child {
        node[probranch] {}
        		child {
        				node[probranch] {}    
        				edge from parent 
           					node[] {}
            				node[]  {}
    				}    
        edge from parent 
            node[] {}
            node[]  {}
    }
    child {
        node[probranch] {}
        		child {
        				node[probranch] {}    
        				edge from parent 
           					node[] {}
            				node[]  {}
    				}    
        edge from parent 
            node[] {}
            node[]  {}
    }
    child {
        node[probranch, label=above: $1$] {}
        		child {
        				node[probranch] {}    
        				edge from parent 
           					node[] {}
            				node[]  {}
    				}
    				child {
        				node[probranch] {}    
        				edge from parent 
           					node[] {}
            				node[]  {}
    				}    
    				child {
        				node[probranch, label=above: $4/3$] {}    
        				edge from parent 
           					node[] {}
            				node[]  {}
    				}          
        edge from parent         
            node[] {}
            node[]  {}
    };
\end{tikzpicture} \\
\end{tabular}
\tikzstyle{level 1}=[level distance=9mm, sibling distance=9mm]
\tikzstyle{level 2}=[level distance=6mm, sibling distance=3mm]
\begin{tabular}{l r}
\begin{tikzpicture}[grow=right, sloped]
\draw  [-, rounded corners, thick] (1.6,.9) -- (1.8,1.05) -- (1.6,1.2);
\draw [-, rounded corners, thick] (1, -.05) -- (1.2, -.45) -- (1, -.85);
\node[probranch] {}
    child {
        node[probranch] {}
        		child {
        				node[probranch] {}    
        				edge from parent 
           					node[] {}
            				node[]  {}
    				}    
        edge from parent 
            node[] {}
            node[]  {}
    }
    child {
        node[probranch] {}
        		child {
        				node[probranch] {}    
        				edge from parent 
           					node[] {}
            				node[]  {}
    				}    
        edge from parent 
            node[] {}
            node[]  {}
    }
    child {
        node[probranch, label=above: $1$] {}
        		child {
        				node[probranch] {}    
        				edge from parent 
           					node[] {}
            				node[]  {}
    				}
    				child {
        				node[probranch] {}    
        				edge from parent 
           					node[] {}
            				node[]  {}
    				}    
    				child {
        				node[probranch, label=above: $4/3$] {}    
        				edge from parent 
           					node[] {}
            				node[]  {}
    				}          
        edge from parent         
            node[] {}
            node[]  {}
    };
\end{tikzpicture} \\
\end{tabular}
\tikzstyle{level 1}=[level distance=6mm, sibling distance=7mm]
\tikzstyle{level 2}=[level distance=6mm, sibling distance=4mm]
\tikzstyle{level 3}=[level distance=6mm, sibling distance=2mm]
\begin{tabular}{l r}
\begin{tikzpicture}[grow=right, sloped]
\draw  [-, rounded corners, thick] (1.9,-.325) -- (2.1,-.45) -- (1.9,-.6);
\node[probranch] {}
    child {
        node[probranch] {}
        		child {
        				node[probranch] {}
        						child {
        								node[probranch] {}    
        						edge from parent 
           							node[] {}
            						node[]  {}
    								}
    								child {
        								node[probranch] {}    
        						edge from parent 
           							node[] {}
            						node[]  {}
    								}      
							child {
        								node[probranch] {}    
        						edge from parent 
           							node[] {}
            						node[]  {}
    								}
        				edge from parent 
           					node[] {}
            				node[]  {}
    				}   
        edge from parent 
            node[] {}
            node[]  {}
    }
    child {
        node[probranch, label=above: $1$] {}
            child {
        				node[probranch] {}
        						child {
        								node[probranch] {}    
        						edge from parent 
           							node[] {}
            						node[]  {}
    								}    
        				edge from parent 
           					node[] {}
            				node[]  {}
    				} 
    				child {
        				node[probranch, label=above: $3/2$] {}
        						child {
        								node[probranch, label=above: $4/3$] {}    
        						edge from parent 
           							node[] {}
            						node[]  {}
    								}    
        				edge from parent 
           					node[] {}
            				node[]  {}
    				}          
        edge from parent         
            node[] {}
            node[]  {}
    };
\end{tikzpicture} \\
\end{tabular}
\tikzstyle{level 1}=[level distance=6mm, sibling distance=7mm]
\tikzstyle{level 2}=[level distance=6mm, sibling distance=4mm]
\tikzstyle{level 3}=[level distance=6mm, sibling distance=2mm]
\begin{tabular}{l r}
\begin{tikzpicture}[grow=right, sloped]
\draw  [-, rounded corners, thick] (1.3,0) -- (1.5,.15) -- (1.3,.3);
\node[probranch] {}
    child {
        node[probranch] {}
        		child {
        				node[probranch] {}
        						child {
        								node[probranch] {}    
        						edge from parent 
           							node[] {}
            						node[]  {}
    								}
    								child {
        								node[probranch] {}    
        						edge from parent 
           							node[] {}
            						node[]  {}
    								}      
        				edge from parent 
           					node[] {}
            				node[]  {}
    				}   
        edge from parent 
            node[] {}
            node[]  {}
    }
    child {
        node[probranch, label=above: $1$] {}
            child {
        				node[probranch] {}
        						child {
        								node[probranch] {}    
        						edge from parent 
           							node[] {}
            						node[]  {}
    								}    
        				edge from parent 
           					node[] {}
            				node[]  {}
    				} 
        		child {
        				node[probranch] {}
        						child {
        								node[probranch] {}    
        						edge from parent 
           							node[] {}
            						node[]  {}
    								}    
        				edge from parent 
           					node[] {}
            				node[]  {}
    				} 
    				child {
        				node[probranch, label=above: $4/3$] {}
        						child {
        								node[probranch, label=above: $c$] {}    
        						edge from parent 
           							node[] {}
            						node[]  {}
    								}    
        				edge from parent 
           					node[] {}
            				node[]  {}
    				}          
        edge from parent         
            node[] {}
            node[]  {}
    };
\end{tikzpicture}\\
\end{tabular}
\tikzstyle{level 1}=[level distance=6mm, sibling distance=7mm]
\tikzstyle{level 2}=[level distance=6mm, sibling distance=4mm]
\tikzstyle{level 3}=[level distance=6mm, sibling distance=2mm]
\begin{tabular}{l r}
\begin{tikzpicture}[grow=right, sloped]
\draw  [-, rounded corners, thick] (1.3,0) -- (1.5,.15) -- (1.3,.3);
\draw [-, rounded corners, thick] (1.9 , -.25) -- (2.1 , -.35) -- (1.9 , -.45);
\node[probranch] {}
    child {
        node[probranch] {}
        		child {
        				node[probranch] {}
        						child {
        								node[probranch] {}    
        						edge from parent 
           							node[] {}
            						node[]  {}
    								}
    								child {
        								node[probranch] {}    
        						edge from parent 
           							node[] {}
            						node[]  {}
    								}      
        				edge from parent 
           					node[] {}
            				node[]  {}
    				}   
        edge from parent 
            node[] {}
            node[]  {}
    }
    child {
        node[probranch, label=above: $1$] {}
            child {
        				node[probranch] {}
        						child {
        								node[probranch] {}    
        						edge from parent 
           							node[] {}
            						node[]  {}
    								}    
        				edge from parent 
           					node[] {}
            				node[]  {}
    				} 
        		child {
        				node[probranch] {}
        						child {
        								node[probranch] {}    
        						edge from parent 
           							node[] {}
            						node[]  {}
    								}    
        				edge from parent 
           					node[] {}
            				node[]  {}
    				} 
    				child {
        				node[probranch, label=above: $4/3$] {}
        						child {
        								node[probranch, label=above: $d$] {}    
        						edge from parent 
           							node[] {}
            						node[]  {}
    								}    
        				edge from parent 
           					node[] {}
            				node[]  {}
    				}          
        edge from parent         
            node[] {}
            node[]  {}
    };
\end{tikzpicture} \\
\end{tabular}
\end{center}
\begin{center}
$a={\left\{\frac{n}{2}\right\}}^{11}_4$, $b={\left\{n\right\}}^5_2$, $c={\left\{\frac{n}{2}\right\}}^{13}_3$, $d={\left\{n\right\}}^6_2$ \\ 
\end{center}
Newton Polygon 17: \\ 
\begin{center}
\begin{tikzpicture}[scale=.25]
\draw[help lines] (0,0) grid (7,7);
\draw(1,5) -- (3,2);
\draw(3,2) -- (5,0);
\draw(0,7)--(1,5);
\draw[fill] (0,7) circle [radius= 0.25];
\draw[fill] (1,6) circle [radius= 0.25];
\draw[fill] (2,5) circle [radius= 0.25];
\draw[fill] (3,4) circle [radius= 0.25];
\draw[fill] (4,3) circle [radius= 0.25];
\draw[fill] (5,2) circle [radius= 0.25];
\draw[fill] (6,1) circle [radius= 0.25];
\draw[fill] (7,0) circle [radius= 0.25];
\draw[fill] (6,0) circle [radius= 0.25];
\draw[fill](5,1) circle [radius = 0.25];
\draw[fill](4,2) circle [radius= 0.25];
\draw[fill](3,3) circle [radius = 0.25];
\draw[fill](2,4) circle[radius = 0.25];
\draw[fill](5,0) circle [radius = 0.25];
\draw[fill](4,1) circle [radius = 0.25];
\draw[fill](1,5) circle [radius = 0.25];
\draw[fill](3,2) circle [radius = 0.25];
\end{tikzpicture}
\end{center}
Families:\\ $y^3(y-x)^2+ax^5y+bx^4y^2+cx^3y^3+dx^2y^4+exy^5+fy^6+gx^7+hx^6y+jx^5y^2+kx^4y^3+lx^3y^4+mx^2y^5+nxy^6+py^7=0$, \\ \\ $y^3(y-x)(y-2x)+ax^5y+bx^4y^2+cx^3y^3+dx^2y^4+exy^5+fy^6+gx^7+hx^6y+jx^5y^2+kx^4y^3+lx^3y^4+mx^2y^5+nxy^6+py^7=0$, \\ \\ $y^3(y^2+x^2)+ax^5y+bx^4y^2+cx^3y^3+dx^2y^4+exy^5+fy^6+gx^7+hx^6y+jx^5y^2+kx^4y^3+lx^3y^4+mx^2y^5+nxy^6+py^7=0$ \\ \\
Diagrams: 
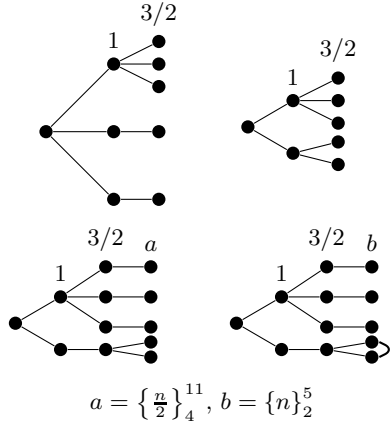
\begin{center}
\tikzstyle{level 1}=[level distance=9mm, sibling distance=9mm]
\tikzstyle{level 2}=[level distance=6mm, sibling distance=3mm]
\begin{tabular}{l r}
\begin{tikzpicture}[grow=right, sloped]
\node[probranch] {}
    child {
        node[probranch] {}
        		child {
        				node[probranch] {}    
        				edge from parent 
           					node[] {}
            				node[]  {}
    				}    
        edge from parent 
            node[] {}
            node[]  {}
    }
    child {
        node[probranch] {}
        		child {
        				node[probranch] {}    
        				edge from parent 
           					node[] {}
            				node[]  {}
    				}    
        edge from parent 
            node[] {}
            node[]  {}
    }
    child {
        node[probranch, label=above: $1$] {}
        		child {
        				node[probranch] {}    
        				edge from parent 
           					node[] {}
            				node[]  {}
    				}
    				child {
        				node[probranch] {}    
        				edge from parent 
           					node[] {}
            				node[]  {}
    				}    
    				child {
        				node[probranch, label=above: $3/2$] {}    
        				edge from parent 
           					node[] {}
            				node[]  {}
    				}          
        edge from parent         
            node[] {}
            node[]  {}
    };
\end{tikzpicture}\\
\end{tabular}
\tikzstyle{level 1}=[level distance=6mm, sibling distance=7mm]
\tikzstyle{level 2}=[level distance=6mm, sibling distance=3mm]
\begin{tabular}{l r}
\begin{tikzpicture}[grow=right, sloped]
\node[probranch] {}
    child {
        node[probranch] {}
        		child {
        				node[probranch] {}    
        				edge from parent 
           					node[] {}
            				node[]  {}
    				}
    				child {
        				node[probranch] {}    
        				edge from parent 
           					node[] {}
            				node[]  {}
    				}     
        edge from parent 
            node[] {}
            node[]  {}
    }
    child {
        node[probranch, label=above: $1$] {}
            child {
        				node[probranch] {}    
        				edge from parent 
           					node[] {}
            				node[]  {}
    				} 
        		child {
        				node[probranch] {}    
        				edge from parent 
           					node[] {}
            				node[]  {}
    				} 
    				child {
        				node[probranch, label=above: $3/2$] {}    
        				edge from parent 
           					node[] {}
            				node[]  {}
    				}          
        edge from parent         
            node[] {}
            node[]  {}
    };
\end{tikzpicture} \\
\end{tabular}
\\
\tikzstyle{level 1}=[level distance=6mm, sibling distance=7mm]
\tikzstyle{level 2}=[level distance=6mm, sibling distance=4mm]
\tikzstyle{level 3}=[level distance=6mm, sibling distance=2mm]
\begin{tabular}{l r}
\begin{tikzpicture}[grow=right, sloped]
\node[probranch] {}
    child {
        node[probranch] {}
        		child {
        				node[probranch] {}
        						child {
        								node[probranch] {}    
        						edge from parent 
           							node[] {}
            						node[]  {}
    								}
    								child {
        								node[probranch] {}    
        						edge from parent 
           							node[] {}
            						node[]  {}
    								}      
        				edge from parent 
           					node[] {}
            				node[]  {}
    				}   
        edge from parent 
            node[] {}
            node[]  {}
    }
    child {
        node[probranch, label=above: $1$] {}
            child {
        				node[probranch] {}
        						child {
        								node[probranch] {}    
        						edge from parent 
           							node[] {}
            						node[]  {}
    								}    
        				edge from parent 
           					node[] {}
            				node[]  {}
    				} 
        		child {
        				node[probranch] {}
        						child {
        								node[probranch] {}    
        						edge from parent 
           							node[] {}
            						node[]  {}
    								}    
        				edge from parent 
           					node[] {}
            				node[]  {}
    				} 
    				child {
        				node[probranch, label=above: $3/2$] {}
        						child {
        								node[probranch, label=above: $a$] {}    
        						edge from parent 
           							node[] {}
            						node[]  {}
    								}    
        				edge from parent 
           					node[] {}
            				node[]  {}
    				}          
        edge from parent         
            node[] {}
            node[]  {}
    };
\end{tikzpicture} \\
\end{tabular}
\tikzstyle{level 1}=[level distance=6mm, sibling distance=7mm]
\tikzstyle{level 2}=[level distance=6mm, sibling distance=4mm]
\tikzstyle{level 3}=[level distance=6mm, sibling distance=2mm]
\begin{tabular}{l r}
\begin{tikzpicture}[grow=right, sloped]
\draw [-, rounded corners, thick] (1.9 , -.25) -- (2.1 , -.35) -- (1.9 , -.45);
\node[probranch] {}
    child {
        node[probranch] {}
        		child {
        				node[probranch] {}
        						child {
        								node[probranch] {}    
        						edge from parent 
           							node[] {}
            						node[]  {}
    								}
    								child {
        								node[probranch] {}    
        						edge from parent 
           							node[] {}
            						node[]  {}
    								}      
        				edge from parent 
           					node[] {}
            				node[]  {}
    				}   
        edge from parent 
            node[] {}
            node[]  {}
    }
    child {
        node[probranch, label=above: $1$] {}
            child {
        				node[probranch] {}
        						child {
        								node[probranch] {}    
        						edge from parent 
           							node[] {}
            						node[]  {}
    								}    
        				edge from parent 
           					node[] {}
            				node[]  {}
    				} 
        		child {
        				node[probranch] {}
        						child {
        								node[probranch] {}    
        						edge from parent 
           							node[] {}
            						node[]  {}
    								}    
        				edge from parent 
           					node[] {}
            				node[]  {}
    				} 
    				child {
        				node[probranch, label=above: $3/2$] {}
        						child {
        								node[probranch, label=above: $b$] {}    
        						edge from parent 
           							node[] {}
            						node[]  {}
    								}    
        				edge from parent 
           					node[] {}
            				node[]  {}
    				}          
        edge from parent         
            node[] {}
            node[]  {}
    };
\end{tikzpicture} \\
\end{tabular}
\end{center}
\begin{center}
$a={\left\{\frac{n}{2}\right\}}^{11}_4$, $b={\left\{n\right\}}^5_2$ \\ 
\end{center}
Tangent Cones: \\ $y^2(y-x)^2(y-2x)$, \\ $y^2(y-x)(y-2x)(y-ax)$ where $a\neq0,1,2$, \\ $y^2(y-x)(y^2+x^2)$ \\ \\
Newton Polygons 18, 19: \\ 
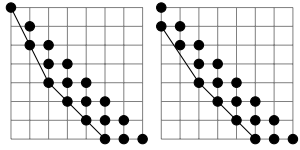
\begin{center}
\begin{tikzpicture}[scale=.25]
\draw[help lines] (0,0) grid (7,7);
\draw(2,3) -- (5,0);
\draw(0,7)--(2,3);
\draw[fill] (0,7) circle [radius= 0.25];
\draw[fill] (1,6) circle [radius= 0.25];
\draw[fill] (2,5) circle [radius= 0.25];
\draw[fill] (3,4) circle [radius= 0.25];
\draw[fill] (4,3) circle [radius= 0.25];
\draw[fill] (5,2) circle [radius= 0.25];
\draw[fill] (6,1) circle [radius= 0.25];
\draw[fill] (7,0) circle [radius= 0.25];
\draw[fill] (6,0) circle [radius= 0.25];
\draw[fill](5,1) circle [radius = 0.25];
\draw[fill](4,2) circle [radius= 0.25];
\draw[fill](3,3) circle [radius = 0.25];
\draw[fill](2,4) circle[radius = 0.25];
\draw[fill](5,0) circle [radius = 0.25];
\draw[fill](4,1) circle [radius = 0.25];
\draw[fill](1,5) circle [radius = 0.25];
\draw[fill](3,2) circle [radius = 0.25];
\draw[fill](2,3) circle [radius = 0.25];
\end{tikzpicture}
\begin{tikzpicture}[scale=.25]
\draw[help lines] (0,0) grid (7,7);
\draw(2,3) -- (5,0);
\draw(0,6)--(2,3);
\draw[fill] (0,7) circle [radius= 0.25];
\draw[fill] (1,6) circle [radius= 0.25];
\draw[fill] (2,5) circle [radius= 0.25];
\draw[fill] (3,4) circle [radius= 0.25];
\draw[fill] (4,3) circle [radius= 0.25];
\draw[fill] (5,2) circle [radius= 0.25];
\draw[fill] (6,1) circle [radius= 0.25];
\draw[fill] (7,0) circle [radius= 0.25];
\draw[fill] (6,0) circle [radius= 0.25];
\draw[fill](5,1) circle [radius = 0.25];
\draw[fill](4,2) circle [radius= 0.25];
\draw[fill](3,3) circle [radius = 0.25];
\draw[fill](2,4) circle[radius = 0.25];
\draw[fill](5,0) circle [radius = 0.25];
\draw[fill](4,1) circle [radius = 0.25];
\draw[fill](1,5) circle [radius = 0.25];
\draw[fill](3,2) circle [radius = 0.25];
\draw[fill](2,3) circle [radius = 0.25];
\draw[fill](0,6) circle [radius = 0.25];
\end{tikzpicture}
\end{center}
Families: \\ $y^2(y-x)^2(y-2x)+ax^5y+bx^4y^2+cx^3y^3+dx^2y^4+exy^5+fy^6+gx^7+hx^6y+jx^5y^2+kx^4y^3+lx^3y^4+mx^2y^5+nxy^6+py^7=0$, \\ \\ $y^2(y-x)(y-2x)(y-ax)+bx^5y+cx^4y^2+dx^3y^3+ex^2y^4+fxy^5+gy^6+hx^7+jx^6y+kx^5y^2+lx^4y^3+mx^3y^4+nx^2y^5+pxy^6+qy^7=0$ where $a\neq0,1,2$, \\ \\ $y^2(y-x)(y^2+x^2)+ax^5y+bx^4y^2+cx^3y^3+dx^2y^4+exy^5+fy^6+gx^7+hx^6y+jx^5y^2+kx^4y^3+lx^3y^4+mx^2y^5+nxy^6+py^7=0$, \\ \\ $y^2(y-x)^2(y-2x)+ax^6+bx^5y+cx^4y^2+dx^3y^3+ex^2y^4+fxy^5+gy^6+hx^7+jx^6y+kx^5y^2+lx^4y^3+mx^3y^4+nx^2y^5+pxy^6+qy^7=0$, \\ \\ $y^2(y-x)(y-2x)(y-ax)+bx^6+cx^5y+dx^4y^2+ex^3y^3+fx^2y^4+gxy^5+hy^6+jx^7+kx^6y+lx^5y^2+mx^4y^3+nx^3y^4+px^2y^5+qxy^6+ry^7=0$ where $a\neq0,1,2$, \\ \\ $y^2(y-x)(y^2+x^2)+ax^6+bx^5y+cx^4y^2+dx^3y^3+ex^2y^4+fxy^5+gy^6+hx^7+jx^6y+kx^5y^2+lx^4y^3+mx^3y^4+nx^2y^5+pxy^6+qy^7=0$ \\ \\
Diagrams: 
\\
\begin{center}
\tikzstyle{level 1}=[level distance=9mm, sibling distance=6mm]
\tikzstyle{level 2}=[level distance=6mm, sibling distance=4mm]
\begin{tabular}{l r}
\begin{tikzpicture}[grow=right, sloped]
\node[probranch] {}
    child {
        node[probranch] {} 
           child {
        				node[probranch] {}    
        				edge from parent 
           					node[] {}
            				node[]  {}
            }     
        edge from parent 
            node[] {}
            node[]  {}
    }
    child {
        node[probranch] {}
        		child {
        				node[probranch] {}    
        				edge from parent 
           					node[] {}
            				node[]  {}
    				}    
        edge from parent 
            node[] {}
            node[]  {}
    }
    child {
        node[probranch] {}
        		child {
        				node[probranch] {}    
        				edge from parent 
           					node[] {}
            				node[]  {}
    				}    
        edge from parent 
            node[] {}
            node[]  {}
    }
    child {
        node[probranch, label=above: $1$] {}
        		child {
        				node[probranch] {}    
        				edge from parent 
           					node[] {}
            				node[]  {}
    				}
    				child {
        				node[probranch, label=above: $a$] {}    
        				edge from parent 
           					node[] {}
            				node[]  {}
    				}          
        edge from parent         
            node[] {}
            node[]  {}
    };
\end{tikzpicture} \\
\end{tabular}
\tikzstyle{level 1}=[level distance=9mm, sibling distance=6mm]
\tikzstyle{level 2}=[level distance=6mm, sibling distance=4mm]
\begin{tabular}{l r}
\begin{tikzpicture}[grow=right, sloped]
\draw  [-, rounded corners, thick] (1.6,.7) -- (1.8,.9) -- (1.6,1.1);
\node[probranch] {}
    child {
        node[probranch] {} 
           child {
        				node[probranch] {}    
        				edge from parent 
           					node[] {}
            				node[]  {}
            }     
        edge from parent 
            node[] {}
            node[]  {}
    }
    child {
        node[probranch] {}
        		child {
        				node[probranch] {}    
        				edge from parent 
           					node[] {}
            				node[]  {}
    				}    
        edge from parent 
            node[] {}
            node[]  {}
    }
    child {
        node[probranch] {}
        		child {
        				node[probranch] {}    
        				edge from parent 
           					node[] {}
            				node[]  {}
    				}    
        edge from parent 
            node[] {}
            node[]  {}
    }
    child {
        node[probranch, label=above: $1$] {}
        		child {
        				node[probranch] {}    
        				edge from parent 
           					node[] {}
            				node[]  {}
    				}
    				child {
        				node[probranch, label=above: $a$] {}    
        				edge from parent 
           					node[] {}
            				node[]  {}
    				}          
        edge from parent         
            node[] {}
            node[]  {}
    };
\end{tikzpicture} \\
\end{tabular}
\tikzstyle{level 1}=[level distance=9mm, sibling distance=6mm]
\tikzstyle{level 2}=[level distance=6mm, sibling distance=4mm]
\begin{tabular}{l r}
\begin{tikzpicture}[grow=right, sloped]
\draw  [-, rounded corners, thick] (1,-.35) -- (1.2,-.6) -- (1,-.85);
\node[probranch] {}
    child {
        node[probranch] {} 
           child {
        				node[probranch] {}    
        				edge from parent 
           					node[] {}
            				node[]  {}
            }     
        edge from parent 
            node[] {}
            node[]  {}
    }
    child {
        node[probranch] {}
        		child {
        				node[probranch] {}    
        				edge from parent 
           					node[] {}
            				node[]  {}
    				}    
        edge from parent 
            node[] {}
            node[]  {}
    }
    child {
        node[probranch] {}
        		child {
        				node[probranch] {}    
        				edge from parent 
           					node[] {}
            				node[]  {}
    				}    
        edge from parent 
            node[] {}
            node[]  {}
    }
    child {
        node[probranch, label=above: $1$] {}
        		child {
        				node[probranch] {}    
        				edge from parent 
           					node[] {}
            				node[]  {}
    				}
    				child {
        				node[probranch, label=above: $a$] {}    
        				edge from parent 
           					node[] {}
            				node[]  {}
    				}          
        edge from parent         
            node[] {}
            node[]  {}
    };
\end{tikzpicture} \\
\end{tabular}
\tikzstyle{level 1}=[level distance=9mm, sibling distance=6mm]
\tikzstyle{level 2}=[level distance=6mm, sibling distance=4mm]
\begin{tabular}{l r}
\begin{tikzpicture}[grow=right, sloped]
\draw  [-, rounded corners, thick] (1.6,.7) -- (1.8,.9) -- (1.6,1.1);
\draw  [-, rounded corners, thick] (1,-.35) -- (1.2,-.6) -- (1,-.85);
\node[probranch] {}
    child {
        node[probranch] {} 
           child {
        				node[probranch] {}    
        				edge from parent 
           					node[] {}
            				node[]  {}
            }     
        edge from parent 
            node[] {}
            node[]  {}
    }
    child {
        node[probranch] {}
        		child {
        				node[probranch] {}    
        				edge from parent 
           					node[] {}
            				node[]  {}
    				}    
        edge from parent 
            node[] {}
            node[]  {}
    }
    child {
        node[probranch] {}
        		child {
        				node[probranch] {}    
        				edge from parent 
           					node[] {}
            				node[]  {}
    				}    
        edge from parent 
            node[] {}
            node[]  {}
    }
    child {
        node[probranch, label=above: $1$] {}
        		child {
        				node[probranch] {}    
        				edge from parent 
           					node[] {}
            				node[]  {}
    				}
    				child {
        				node[probranch, label=above: $a$] {}    
        				edge from parent 
           					node[] {}
            				node[]  {}
    				}          
        edge from parent         
            node[] {}
            node[]  {}
    };
\end{tikzpicture} \\
\end{tabular}
\\
$a={\left\{\frac{n}{2}\right\}}^{13}_3$ \hspace{7 mm} $a={\left\{n\right\}}^{6}_2$ \hspace{7 mm} $a={\left\{\frac{n}{2}\right\}}^{13}_4$ \hspace{7 mm} $a={\left\{n\right\}}^{6}_2$
\\
\tikzstyle{level 1}=[level distance=4mm, sibling distance=7mm]
\tikzstyle{level 2}=[level distance=6mm, sibling distance=4mm]
\tikzstyle{level 3}=[level distance=6mm, sibling distance=3mm]
\begin{tabular}{l r}
\begin{tikzpicture}[grow=right, sloped]
\node[probranch] {}
    child {
        node[probranch] {}
        		child {
        				node[probranch] {}
        						child {
        								node[probranch] {}    
        						edge from parent 
           							node[] {}
            						node[]  {}
    								}   
        				edge from parent 
           					node[] {}
            				node[]  {}
    				}   
        edge from parent 
            node[] {}
            node[]  {}
    }
    child {
        node[probranch] {}
        		child {
        				node[probranch] {}
        						child {
        								node[probranch] {}    
        						edge from parent 
           							node[] {}
            						node[]  {}
    								}
    								child {
        								node[probranch] {}    
        						edge from parent 
           							node[] {}
            						node[]  {}
    								}      
        				edge from parent 
           					node[] {}
            				node[]  {}
    				}   
        edge from parent 
            node[] {}
            node[]  {}
    }
    child {
        node[probranch, label=above: $1$] {}
        		child {
        				node[probranch] {}
        						child {
        								node[probranch] {}    
        						edge from parent 
           							node[] {}
            						node[]  {}
    								}    
        				edge from parent 
           					node[] {}
            				node[]  {}
    				} 
    				child {
        				node[probranch, label=above: $3/2$] {}
        						child {
        								node[probranch, label=above: $a$] {}    
        						edge from parent 
           							node[] {}
            						node[]  {}
    								}    
        				edge from parent 
           					node[] {}
            				node[]  {}
    				}          
        edge from parent         
            node[] {}
            node[]  {}
    };
\end{tikzpicture} \\
\end{tabular}
\tikzstyle{level 1}=[level distance=4mm, sibling distance=7mm]
\tikzstyle{level 2}=[level distance=6mm, sibling distance=4mm]
\tikzstyle{level 3}=[level distance=6mm, sibling distance=3mm]
\begin{tabular}{l r}
\begin{tikzpicture}[grow=right, sloped]
\draw  [-, rounded corners, thick] (1.7,.15) -- (1.9,0) -- (1.7,-.15);
\node[probranch] {}
    child {
        node[probranch] {}
        		child {
        				node[probranch] {}
        						child {
        								node[probranch] {}    
        						edge from parent 
           							node[] {}
            						node[]  {}
    								}   
        				edge from parent 
           					node[] {}
            				node[]  {}
    				}   
        edge from parent 
            node[] {}
            node[]  {}
    }
    child {
        node[probranch] {}
        		child {
        				node[probranch] {}
        						child {
        								node[probranch] {}    
        						edge from parent 
           							node[] {}
            						node[]  {}
    								}
    								child {
        								node[probranch] {}    
        						edge from parent 
           							node[] {}
            						node[]  {}
    								}      
        				edge from parent 
           					node[] {}
            				node[]  {}
    				}   
        edge from parent 
            node[] {}
            node[]  {}
    }
    child {
        node[probranch, label=above: $1$] {}
        		child {
        				node[probranch] {}
        						child {
        								node[probranch] {}    
        						edge from parent 
           							node[] {}
            						node[]  {}
    								}    
        				edge from parent 
           					node[] {}
            				node[]  {}
    				} 
    				child {
        				node[probranch, label=above: $3/2$] {}
        						child {
        								node[probranch, label=above: $a$] {}    
        						edge from parent 
           							node[] {}
            						node[]  {}
    								}    
        				edge from parent 
           					node[] {}
            				node[]  {}
    				}          
        edge from parent         
            node[] {}
            node[]  {}
    };
\end{tikzpicture}\\
\end{tabular}
\tikzstyle{level 1}=[level distance=4mm, sibling distance=7mm]
\tikzstyle{level 2}=[level distance=6mm, sibling distance=4mm]
\tikzstyle{level 3}=[level distance=6mm, sibling distance=3mm]
\begin{tabular}{l r}
\begin{tikzpicture}[grow=right, sloped]
\node[probranch] {}
    child {
        node[probranch] {}
        		child {
        				node[probranch] {}
        						child {
        								node[probranch] {}    
        						edge from parent 
           							node[] {}
            						node[]  {}
    								}   
        				edge from parent 
           					node[] {}
            				node[]  {}
    				}   
        edge from parent 
            node[] {}
            node[]  {}
    }
    child {
        node[probranch] {}
        		child {
        				node[probranch] {}
        						child {
        								node[probranch] {}    
        						edge from parent 
           							node[] {}
            						node[]  {}
    								}
    								child {
        								node[probranch] {}    
        						edge from parent 
           							node[] {}
            						node[]  {}
    								}      
        				edge from parent 
           					node[] {}
            				node[]  {}
    				}   
        edge from parent 
            node[] {}
            node[]  {}
    }
    child {
        node[probranch, label=above: $1$] {}
        		child {
        				node[probranch] {}
        						child {
        								node[probranch] {}    
        						edge from parent 
           							node[] {}
            						node[]  {}
    								}    
        				edge from parent 
           					node[] {}
            				node[]  {}
    				} 
    				child {
        				node[probranch, label=above: $2$] {}
        						child {
        								node[probranch, label=above: $a$] {}    
        						edge from parent 
           							node[] {}
            						node[]  {}
    								}    
        				edge from parent 
           					node[] {}
            				node[]  {}
    				}          
        edge from parent         
            node[] {}
            node[]  {}
    };
\end{tikzpicture} \\
\end{tabular}\\
$a={\left\{\frac{n}{2}\right\}}^{13}_4$ \hspace{14 mm} $a={\left\{n\right\}}^{6}_2$ \hspace{14 mm} $a={\left\{\frac{n}{2}\right\}}^{11}_5$
\\
\tikzstyle{level 1}=[level distance=4mm, sibling distance=7mm]
\tikzstyle{level 2}=[level distance=6mm, sibling distance=4mm]
\tikzstyle{level 3}=[level distance=6mm, sibling distance=3mm]
%
\tikzstyle{level 1}=[level distance=4mm, sibling distance=7mm]
\tikzstyle{level 2}=[level distance=6mm, sibling distance=4mm]
\tikzstyle{level 3}=[level distance=6mm, sibling distance=3mm]
%
\\
$a={\left\{\frac{n}{2}\right\}}^{11}_5$ \hspace{14 mm} $a={\left\{n\right\}}^{5}_3$ \hspace{14 mm}  $a={\left\{n\right\}}^{5}_3$
\\
\tikzstyle{level 1}=[level distance=4mm, sibling distance=7mm]
\tikzstyle{level 2}=[level distance=6mm, sibling distance=4mm]
\tikzstyle{level 3}=[level distance=6mm, sibling distance=3mm]
%
\\
$a={\left\{n\right\}}^{5}_3$ \hspace{16 mm} $a={\left\{\frac{n}{2}\right\}}^{11}_6$ \hspace{14 mm} $a={\left\{n\right\}}^{5}_3$
\\
\tikzstyle{level 1}=[level distance=4mm, sibling distance=7mm]
\tikzstyle{level 2}=[level distance=6mm, sibling distance=4mm]
\tikzstyle{level 3}=[level distance=6mm, sibling distance=3mm]
\begin{tabular}{l r}
\begin{tikzpicture}[grow=right, sloped]
\node[probranch] {}
    child {
        node[probranch] {}
        		child {
        				node[probranch] {}
        						child {
        								node[probranch] {}    
        						edge from parent 
           							node[] {}
            						node[]  {}
    								}   
        				edge from parent 
           					node[] {}
            				node[]  {}
    				}   
        edge from parent 
            node[] {}
            node[]  {}
    }
    child {
        node[probranch] {}
        		child {
        				node[probranch] {}
        						child {
        								node[probranch] {}    
        						edge from parent 
           							node[] {}
            						node[]  {}
    								}
    								child {
        								node[probranch] {}    
        						edge from parent 
           							node[] {}
            						node[]  {}
    								}      
        				edge from parent 
           					node[] {}
            				node[]  {}
    				}   
        edge from parent 
            node[] {}
            node[]  {}
    }
    child {
        node[probranch, label=above: $1$] {}
        		child {
        				node[probranch] {}
        						child {
        								node[probranch] {}    
        						edge from parent 
           							node[] {}
            						node[]  {}
    								}    
        				edge from parent 
           					node[] {}
            				node[]  {}
    				} 
    				child {
        				node[probranch, label=above: $3$] {}
        						child {
        								node[probranch, label=above: $a$] {}    
        						edge from parent 
           							node[] {}
            						node[]  {}
    								}    
        				edge from parent 
           					node[] {}
            				node[]  {}
    				}          
        edge from parent         
            node[] {}
            node[]  {}
    };
\end{tikzpicture} \\
\end{tabular}%
\\
$a={\left\{\frac{n}{2}\right\}}^{9}_7$ \hspace{16 mm} $a={\left\{\frac{n}{2}\right\}}^{9}_7$ \hspace{14 mm} $a=4$
\\
\tikzstyle{level 1}=[level distance=4mm, sibling distance=7mm]
\tikzstyle{level 2}=[level distance=6mm, sibling distance=4mm]
\tikzstyle{level 3}=[level distance=6mm, sibling distance=3mm]
%
\\
$a=4$ \hspace{22 mm} $a={\left\{\frac{n}{2}\right\}}^{9}_8$  \hspace{18 mm} $a={\left\{\frac{n}{2}\right\}}^{9}_8$
\\
\tikzstyle{level 1}=[level distance=4mm, sibling distance=7mm]
\tikzstyle{level 2}=[level distance=6mm, sibling distance=4mm]
\tikzstyle{level 3}=[level distance=6mm, sibling distance=3mm]
%

\end{center}
Families: \\ $y(y-x)(y+x)(y-ax)(y-bx)+cx^6+dx^5y+ex^4y^2+fx^3y^3+gx^2y^4+hxy^5+jy^6+kx^7+lx^6y+mx^5y^2+nx^4y^3+px^3y^4+qx^2y^5+rxy^6+sy^7=0$ where $a, b\neq0,1,-1$ and $a\neq b$, \\ \\ $y(y^2+x^2)(y-ax)(y-bx)+cx^6+dx^5y+ex^4y^2+fx^3y^3+gx^2y^4+hxy^5+jy^6+kx^7+lx^6y+mx^5y^2+nx^4y^3+px^3y^4+qx^2y^5+rxy^6+sy^7=0$ where $a, b\neq0$ and $a\neq b$, \\ \\ $y(y^2+x^2)(y^2+ax^2)+bx^6+cx^5y+dx^4y^2+ex^3y^3+fx^2y^4+gxy^5+hy^6+jx^7+kx^6y+lx^5y^2+mx^4y^3+nx^3y^4+px^2y^5+qxy^6+ry^7=0 $ where $a>0$ and $a\neq1$, \\ 
\\$y(y^2+x^2)^2+ax^6+bx^5y+cx^4y^2+dx^3y^3+ex^2y^4+fxy^5+gy^6+hx^7+jx^6y+kx^5y^2+lx^4y^3+mx^3y^4+nx^2y^5+pxy^6+qy^7=0$, 
\\ \\ $y(y-x)(y+x)(y-ax)(y-bx)+cx^5y+dx^4y^2+ex^3y^3+fx^2y^4+gxy^5+hy^6+jx^7+kx^6y+lx^5y^2+mx^4y^3+nx^3y^4+px^2y^5+qxy^6+ry^7=0$ where $a, b\neq0,1,-1$ and $a\neq b$, \\ \\$y(y^2+x^2)(y-ax)(y-bx)+cx^5y+dx^4y^2+ex^3y^3+fx^2y^4+gxy^5+hy^6+jx^7+kx^6y+lx^5y^2+mx^4y^3+nx^3y^4+px^2y^5+qxy^6+ry^7=0$ where $a, b\neq0$ and $a\neq b$,
\\ \\ $y(y^2+x^2)(y^2+ax^2)+bx^5y+cx^4y^2+dx^3y^3+ex^2y^4+fxy^5+gy^6+hx^7+jx^6y+kx^5y^2+lx^4y^3+mx^3y^4+nx^2y^5+pxy^6+qy^7=0 $  where $a>0$ and $a\neq1$, \\ 
\\$y(y^2+x^2)^2+ax^5y+bx^4y^2+cx^3y^3+dx^2y^4+exy^5+fy^6+gx^7+hx^6y+jx^5y^2+kx^4y^3+lx^3y^4+mx^2y^5+nxy^6+py^7=0$ 
\\ \\
Diagrams: 
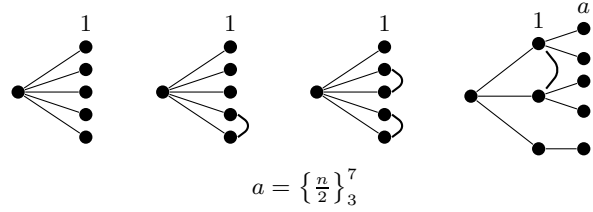
\begin{center}
\tikzstyle{level 1}=[level distance=9mm, sibling distance=3mm]
\tikzstyle{level 2}=[level distance=6mm, sibling distance=6mm]

\tikzstyle{level 1}=[level distance=9mm, sibling distance=3mm]
\tikzstyle{level 2}=[level distance=6mm, sibling distance=6mm]
%
\tikzstyle{level 1}=[level distance=9mm, sibling distance=3mm]
\tikzstyle{level 2}=[level distance=6mm, sibling distance=6mm]
%
\tikzstyle{level 1}=[level distance=9mm, sibling distance=7mm]
\tikzstyle{level 2}=[level distance=6mm, sibling distance=4mm]
%
\\

$a={\left\{\frac{n}{2}\right\}}^{7}_3$
\end{center}

\section{Real Singular Points of Multiplicity 4}

\noindent Tangent Cone: $y^4$ \\ 
\\
Newton Polygon 22: \\ 
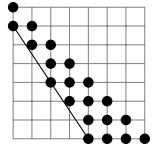
\begin{center}
%
\end{center}
Family: \\ $y^4+ax^2y^3+bxy^4+cy^5+dx^4y^2+ex^3y^3+fx^2y^4+gxy^5+hy^6+jx^7+kx^6y+lx^5y^2+mx^4y^3+nx^3y^4+px^2y^5+qxy^6+ry^7=0$, \\ \\
Diagram: 
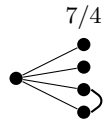
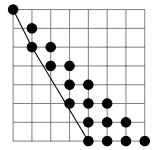
\begin{center}
\tikzstyle{level 1}=[level distance=9mm, sibling distance=3mm]
\tikzstyle{level 2}=[level distance=6mm, sibling distance=6mm]
%
\end{center}
Family: \\ $y^4+ax^3y^2+bx^2y^3+cxy^4+dy^5+ex^6+fx^5y+gx^4y^2+hx^3y^3+jx^2y^4+kxy^5+ly^6+mx^7+nx^6y+px^5y^2+qx^4y^3+rx^3y^4+sx^2y^5+txy^6+uy^7=0$ \\ \\
Diagrams: 
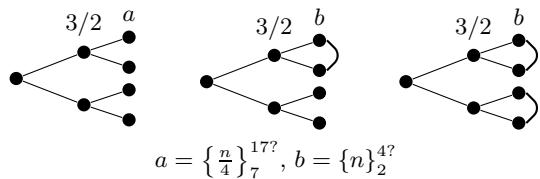
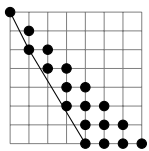
\begin{center}
\tikzstyle{level 1}=[level distance=9mm, sibling distance=7mm]
\tikzstyle{level 2}=[level distance=6mm, sibling distance=4mm]
%
\end{center}
Family: \\ $y^4+ax^2y^3+bxy^4+cy^5+dx^5y+ex^4y^2+fx^3y^3+gx^2y^4+hxy^5+jy^6+kx^7+lx^6y+mx^5y^2+nx^4y^3+px^3y^4+qx^2y^5+rxy^6+sy^7=0$ \\ \\
Diagrams: 
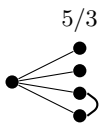
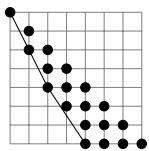
\begin{center}
\tikzstyle{level 1}=[level distance=9mm, sibling distance=3mm]
\tikzstyle{level 2}=[level distance=6mm, sibling distance=6mm]
%
\end{center}
Family: \\ $y^4+ax^3y^2+bx^2y^3+cxy^4+dy^5+ex^5y+fx^4y^2+gx^3y^3+hx^2y^4+jxy^5+ky^6+lx^7+mx^6y+nx^5y^2+px^4y^3+qx^3y^4+rx^2y^5+sxy^6+ty^7=0$ \\ \\
Diagrams: 
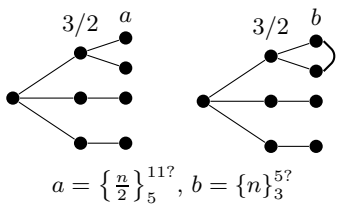
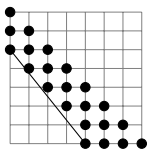
\begin{center}
\tikzstyle{level 1}=[level distance=9mm, sibling distance=6mm]
\tikzstyle{level 2}=[level distance=6mm, sibling distance=4mm]
%
\end{center}
Family: \\ $y^4+ax^5+bx^4y+cx^3y^2+dx^2y^3+exy^4+fy^5+gx^6+hx^5y+jx^4y^2+kx^3y^3+lx^2y^4+mxy^5+ny^6+px^7+qx^6y+rx^5y^2+sx^4y^3+tx^3y^4+ux^2y^5+vxy^6+wy^7=0$ \\ \\
Diagram: 
\begin{center}
\tikzstyle{level 1}=[level distance=9mm, sibling distance=3mm]
\tikzstyle{level 2}=[level distance=6mm, sibling distance=6mm]
%

\end{center}
Families: \\ $y^4+ax^4y+bx^3y^2+cx^2y^3+dxy^4+ey^5+fx^5y+gx^4y^2+hx^3y^3+jx^2y^4+kxy^5+ly^6+mx^7+nx^6y+px^5y^2+qx^4y^3+rx^3y^4+sx^2y^5+txy^6+uy^7=0$, \\ \\ $y^4+ax^4y+bx^3y^2+cx^2y^3+dxy^4+ey^5+fx^6+gx^5y+hx^4y^2+jx^3y^3+kx^2y^4+lxy^5+my^6+nx^7+px^6y+qx^5y^2+rx^4y^3+sx^3y^4+tx^2y^5+uxy^6+vy^7=0$ \\ \\
Diagram: 
\begin{center}
\tikzstyle{level 1}=[level distance=9mm, sibling distance=3mm]
\tikzstyle{level 2}=[level distance=6mm, sibling distance=6mm]
%
\end{center}
Families: \\ $y^3(y-x)+ax^4y+bx^3y^2+cx^2y^3+dxy^4+ey^5+fx^5y+gx^4y^2+hx^3y^3+jx^2y^4+kxy^5+ly^6+mx^7+nx^6y+px^5y^2+qx^4y^3+rx^3y^4+sx^2y^5+txy^6+uy^7=0$, \\ \\ $y^3(y-x)+ax^4y+bx^3y^2+cx^2y^3+dxy^4+ey^5+fx^6+gx^5y+hx^4y^2+jx^3y^3+kx^2y^4+lxy^5+my^6+nx^7+px^6y+qx^5y^2+rx^4y^3+sx^3y^4+tx^2y^5+uxy^6+vy^7=0$ \\ \\
Diagram: 
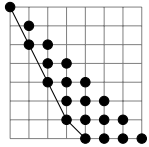
\begin{center}
\tikzstyle{level 1}=[level distance=9mm, sibling distance=9mm]
\tikzstyle{level 2}=[level distance=6mm, sibling distance=3mm]
%
\end{center}
Family: \\ $y^3(y-x)+ax^3y^2+bx^2y^3+cxy^4+dy^5+ex^5y+fx^4y^2+gx^3y^3+hx^2y^4+jxy^5+ky^6+lx^7+mx^6y+nx^5y^2+px^4y^3+qx^3y^4+rx^2y^5+sxy^6+ty^7=0$ \\ \\
Diagrams: \\
\begin{center}
\tikzstyle{level 1}=[level distance=9mm, sibling distance=9mm]
\tikzstyle{level 2}=[level distance=6mm, sibling distance=3mm]
%
\tikzstyle{level 1}=[level distance=6mm, sibling distance=7mm]
\tikzstyle{level 2}=[level distance=6mm, sibling distance=4mm]
\tikzstyle{level 3}=[level distance=6mm, sibling distance=3mm]
%
\tikzstyle{level 1}=[level distance=6mm, sibling distance=7mm]
\tikzstyle{level 2}=[level distance=6mm, sibling distance=4mm]
\tikzstyle{level 3}=[level distance=6mm, sibling distance=3mm]
%
\tikzstyle{level 1}=[level distance=6mm, sibling distance=7mm]
\tikzstyle{level 2}=[level distance=6mm, sibling distance=4mm]
\tikzstyle{level 3}=[level distance=6mm, sibling distance=3mm]
%
\end{center}
Family: \\  $y^3(y-x)+ax^3y^2+bx^2y^3+cxy^4+dy^5+ex^6+fx^5y+gx^4y^2+hx^3y^3+jx^2y^4+kxy^5+ly^6+mx^7+nx^6y+px^5y^2+qx^4y^3+rx^3y^4+sx^2y^5+txy^6+uy^7=0$ \\ \\
Diagram: 
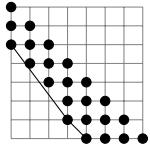
\begin{center}
\tikzstyle{level 1}=[level distance=9mm, sibling distance=9mm]
\tikzstyle{level 2}=[level distance=6mm, sibling distance=3mm]
%
\end{center}
Family: \\  $y^3(y-x)+ax^5+bx^4y+cx^3y^2+dx^2y^3+exy^4+fy^5+gx^6+hx^5y+jx^4y^2+kx^3y^3+lx^2y^4+mxy^5+ny^6+px^7+qx^6y+rx^5y^2+sx^4y^3+tx^3y^4+ux^2y^5+vxy^6+wy^7=0$ \\ \\
Diagram: 
\begin{center}
\tikzstyle{level 1}=[level distance=9mm, sibling distance=9mm]
\tikzstyle{level 2}=[level distance=6mm, sibling distance=3mm]
%
\end{center}
Families: \\ $y^2(x-y)(x+y)+ax^3y^2+bx^2y^3+cxy^4+dy^5+ex^5y+fx^4y^2+gx^3y^3+hx^2y^4+jxy^5+ky^6+lx^7+mx^6y+nx^5y^2+px^4y^3+qx^3y^4+rx^2y^5+sxy^6+ty^7=0$, \\ \\ $y^2(x^2+y^2)+ax^3y^2+bx^2y^3+cxy^4+dy^5+ex^5y+fx^4y^2+gx^3y^3+hx^2y^4+jxy^5+ky^6+lx^7+mx^6y+nx^5y^2+px^4y^3+qx^3y^4+rx^2y^5+sxy^6+ty^7=0$, \\ \\ $y^2(x-ay)^2+bx^3y^2+cx^2y^3+dxy^4+ey^5+fx^5y+gx^4y^2+hx^3y^3+jx^2y^4+kxy^5+ly^6+mx^7+nx^6y+px^5y^2+qx^4y^3+rx^3y^4+sx^2y^5+txy^6+uy^7=0$ where $a\neq0$ \\ \\
Diagrams: \\
No new diagrams. Families are transformed to have Newton polygon 31. \\ 
\\
Newton Polygon 35: \\ 
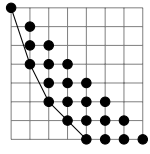
\begin{center}
\begin{tikzpicture}[scale=.25]
\draw[help lines] (0,0) grid (7,7);
\draw(0,7)--(1,4);
\draw(1,4)--(2,2);
\draw(2,2)--(4,0);
\draw[fill] (0,7) circle [radius= 0.25];
\draw[fill] (1,6) circle [radius= 0.25];
\draw[fill] (2,5) circle [radius= 0.25];
\draw[fill] (3,4) circle [radius= 0.25];
\draw[fill] (4,3) circle [radius= 0.25];
\draw[fill] (5,2) circle [radius= 0.25];
\draw[fill] (6,1) circle [radius= 0.25];
\draw[fill] (7,0) circle [radius= 0.25];
\draw[fill] (6,0) circle [radius= 0.25];
\draw[fill](5,1) circle [radius = 0.25];
\draw[fill](4,2) circle [radius= 0.25];
\draw[fill](3,3) circle [radius = 0.25];
\draw[fill](2,4) circle[radius = 0.25];
\draw[fill](5,0) circle [radius = 0.25];
\draw[fill](4,1) circle [radius = 0.25];
\draw[fill](3,2) circle [radius = 0.25];
\draw[fill](4,0) circle [radius = 0.25];
\draw[fill](1,5) circle [radius = 0.25];
\draw[fill](2,3) circle [radius = 0.25];
\draw[fill](3,1) circle [radius = 0.25];
\draw[fill](2,2) circle [radius = 0.25];
\draw[fill](1,4) circle [radius = 0.25];
\end{tikzpicture}
\end{center}
Families: \\ $y^2(x-y)(x+y)+ax^4y+bx^3y^2+cx^2y^3+dxy^4+ey^5+fx^5y+gx^4y^2+hx^3y^3+jx^2y^4+kxy^5+ly^6+mx^7+nx^6y+px^5y^2+qx^4y^3+rx^3y^4+sx^2y^5+txy^6+uy^7=0$, \\ \\ $y^2(x^2+y^2)+ax^4y+bx^3y^2+cx^2y^3+dxy^4+ey^5+fx^5y+gx^4y^2+hx^3y^3+jx^2y^4+kxy^5+ly^6+mx^7+nx^6y+px^5y^2+qx^4y^3+rx^3y^4+sx^2y^5+txy^6+uy^7=0$, \\  \\ $y^2(ax-ay)^2+bx^4y+cx^3y^2+dx^2y^3+exy^4+fy^5+gx^5y+hx^4y^2+jx^3y^3+kx^2y^4+lxy^5+my^6+nx^7+px^6y+qx^5y^2+rx^4y^3+sx^3y^4+tx^2y^5+uxy^6+vy^7=0$ where $a\neq0$ \\ \\
Diagrams: \\
No new diagrams. Families are transformed to have Newton polygon 31. \\ 
\\
Newton Polygon 36: \\ 
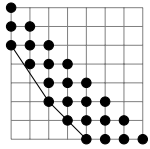
\begin{center}
\begin{tikzpicture}[scale=.25]
\draw[help lines] (0,0) grid (7,7);
\draw(0,5)--(2,2);
\draw(2,2)--(4,0);
\draw[fill] (0,7) circle [radius= 0.25];
\draw[fill] (1,6) circle [radius= 0.25];
\draw[fill] (2,5) circle [radius= 0.25];
\draw[fill] (3,4) circle [radius= 0.25];
\draw[fill] (4,3) circle [radius= 0.25];
\draw[fill] (5,2) circle [radius= 0.25];
\draw[fill] (6,1) circle [radius= 0.25];
\draw[fill] (7,0) circle [radius= 0.25];
\draw[fill] (6,0) circle [radius= 0.25];
\draw[fill](5,1) circle [radius = 0.25];
\draw[fill](4,2) circle [radius= 0.25];
\draw[fill](3,3) circle [radius = 0.25];
\draw[fill](2,4) circle[radius = 0.25];
\draw[fill](5,0) circle [radius = 0.25];
\draw[fill](4,1) circle [radius = 0.25];
\draw[fill](3,2) circle [radius = 0.25];
\draw[fill](4,0) circle [radius = 0.25];
\draw[fill](1,5) circle [radius = 0.25];
\draw[fill](2,3) circle [radius = 0.25];
\draw[fill](3,1) circle [radius = 0.25];
\draw[fill](2,2) circle [radius = 0.25];
\draw[fill](1,4) circle [radius = 0.25];
\draw[fill](0,5) circle [radius = 0.25];
\draw[fill](0,6) circle [radius = 0.25];
\end{tikzpicture}
\end{center}
Families: \\ $y^2(x-y)(x+y)+ax^5+bx^4y+cx^3y^2+dx^2y^3+exy^4+fy^5+gx^6+hx^5y+jx^4y^2+kx^3y^3+lx^2y^4+mxy^5+ny^6+px^7+qx^6y+rx^5y^2+sx^4y^3+tx^3y^4+ux^2y^5+vxy^6+wy^7=0$, \\ \\ $y^2(x^2+y^2)+ax^5+bx^4y+cx^3y^2+dx^2y^3+exy^4+fy^5+gx^6+hx^5y+jx^4y^2+kx^3y^3+lx^2y^4+mxy^5+ny^6+px^7+qx^6y+rx^5y^2+sx^4y^3+tx^3y^4+ux^2y^5+vxy^6+wy^7=0$, \\ \\ $y^2(x-ay)^2+bx^5+cx^4y+dx^3y^2+ex^2y^3+fxy^4+gy^5+hx^6+jx^5y+kx^4y^2+lx^3y^3+mx^2y^4+nxy^5+py^6+qx^7+rx^6y+sx^5y^2+tx^4y^3+ux^3y^4+vx^2y^5+wxy^6+xy^7=0$ where $a\neq0$ \\ \\
Diagrams: 
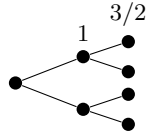
\begin{center}
\tikzstyle{level 1}=[level distance=9mm, sibling distance=7mm]
\tikzstyle{level 2}=[level distance=6mm, sibling distance=4mm]
\begin{tabular}{l r}
\begin{tikzpicture}[grow=right, sloped]
\node[probranch] {}
    child {
        node[probranch] {}
        		child {
        				node[probranch] {}    
        				edge from parent 
           					node[] {}
            				node[]  {}
    				}
    				child {
        				node[probranch] {}    
        				edge from parent 
           					node[] {}
            				node[]  {}
    				}     
        edge from parent 
            node[] {}
            node[]  {}
    }
    child {
        node[probranch, label=above: $1$] {}
        		child {
        				node[probranch] {}    
        				edge from parent 
           					node[] {}
            				node[]  {}
    				} 
    				child {
        				node[probranch, label=above: $3/2$] {}    
        				edge from parent 
           					node[] {}
            				node[]  {}
    				}          
        edge from parent         
            node[] {}
            node[]  {}
    };
\end{tikzpicture} \\
\end{tabular}
\end{center}
Only one new diagram. Families are transformed to have Newton polygon 31. \\ 
\\
Newton Polygon 37: \\ 
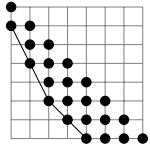
\begin{center}
\begin{tikzpicture}[scale=.25]
\draw[help lines] (0,0) grid (7,7);
\draw(0,6)--(2,2);
\draw(2,2)--(4,0);
\draw[fill] (0,7) circle [radius= 0.25];
\draw[fill] (1,6) circle [radius= 0.25];
\draw[fill] (2,5) circle [radius= 0.25];
\draw[fill] (3,4) circle [radius= 0.25];
\draw[fill] (4,3) circle [radius= 0.25];
\draw[fill] (5,2) circle [radius= 0.25];
\draw[fill] (6,1) circle [radius= 0.25];
\draw[fill] (7,0) circle [radius= 0.25];
\draw[fill] (6,0) circle [radius= 0.25];
\draw[fill](5,1) circle [radius = 0.25];
\draw[fill](4,2) circle [radius= 0.25];
\draw[fill](3,3) circle [radius = 0.25];
\draw[fill](2,4) circle[radius = 0.25];
\draw[fill](5,0) circle [radius = 0.25];
\draw[fill](4,1) circle [radius = 0.25];
\draw[fill](3,2) circle [radius = 0.25];
\draw[fill](4,0) circle [radius = 0.25];
\draw[fill](1,5) circle [radius = 0.25];
\draw[fill](2,3) circle [radius = 0.25];
\draw[fill](3,1) circle [radius = 0.25];
\draw[fill](2,2) circle [radius = 0.25];
\draw[fill](1,4) circle [radius = 0.25];
\draw[fill](0,6) circle [radius = 0.25];
\end{tikzpicture}
\end{center}
Families: \\ $y^2(x-y)(x+y)+ax^4y+bx^3y^2+cx^2y^3+dxy^4+ey^5+fx^6+gx^5y+hx^4y^2+jx^3y^3+kx^2y^4+lxy^5+my^6+nx^7+px^6y+qx^5y^2+rx^4y^3+sx^3y^4+tx^2y^5+uxy^6+vy^7=0$, \\ \\ $y^2(x^2+y^2)+ax^4y+bx^3y^2+cx^2y^3+dxy^4+ey^5+fx^6+gx^5y+hx^4y^2+jx^3y^3+kx^2y^4+lxy^5+my^6+nx^7+px^6y+qx^5y^2+rx^4y^3+sx^3y^4+tx^2y^5+uxy^6+vy^7=0$, \\ \\ $y^2(x-ay)^2+bx^4y+cx^3y^2+dx^2y^3+exy^4+fy^5+gx^6+hx^5y+jx^4y^2+kx^3y^3+lx^2y^4+mxy^5+ny^6+px^7+qx^6y+rx^5y^2+sx^4y^3+tx^3y^4+ux^2y^5+vxy^6+wy^7=0$ where $a\neq0$ \\ \\
Diagrams: \\
\begin{center}
\tikzstyle{level 1}=[level distance=9mm, sibling distance=6mm]
\tikzstyle{level 2}=[level distance=6mm, sibling distance=4mm]
\begin{tabular}{l r}
\begin{tikzpicture}[grow=right, sloped]
\node[probranch] {}
    child {
        node[probranch] {} 
           child {
        				node[probranch] {}    
        				edge from parent 
           					node[] {}
            				node[]  {}
            }     
        edge from parent 
            node[] {}
            node[]  {}
    }
    child {
        node[probranch] {}
        		child {
        				node[probranch] {}    
        				edge from parent 
           					node[] {}
            				node[]  {}
    				}    
        edge from parent 
            node[] {}
            node[]  {}
    }
    child {
        node[probranch, label=above: $1$] {}
        		child {
        				node[probranch] {}    
        				edge from parent 
           					node[] {}
            				node[]  {}
    				}
    				child {
        				node[probranch, label=above: $a$] {}    
        				edge from parent 
           					node[] {}
            				node[]  {}
    				}          
        edge from parent         
            node[] {}
            node[]  {}
    };
\end{tikzpicture} \\
\end{tabular}
\tikzstyle{level 1}=[level distance=9mm, sibling distance=6mm]
\tikzstyle{level 2}=[level distance=6mm, sibling distance=4mm]
\begin{tabular}{l r}
\begin{tikzpicture}[grow=right, sloped]
\draw  [-, rounded corners, thick] (1,-.05) -- (1.2,-.3) -- (1,-.55);
\node[probranch] {}
    child {
        node[probranch] {} 
           child {
        				node[probranch] {}    
        				edge from parent 
           					node[] {}
            				node[]  {}
            }     
        edge from parent 
            node[] {}
            node[]  {}
    }
    child {
        node[probranch] {}
        		child {
        				node[probranch] {}    
        				edge from parent 
           					node[] {}
            				node[]  {}
    				}    
        edge from parent 
            node[] {}
            node[]  {}
    }
    child {
        node[probranch, label=above: $1$] {}
        		child {
        				node[probranch] {}    
        				edge from parent 
           					node[] {}
            				node[]  {}
    				}
    				child {
        				node[probranch, label=above: $b$] {}    
        				edge from parent 
           					node[] {}
            				node[]  {}
    				}          
        edge from parent         
            node[] {}
            node[]  {}
    };
\end{tikzpicture} \\
\end{tabular}
\tikzstyle{level 1}=[level distance=9mm, sibling distance=6mm]
\tikzstyle{level 2}=[level distance=6mm, sibling distance=4mm]
\begin{tabular}{l r}
\begin{tikzpicture}[grow=right, sloped]
\draw  [-, rounded corners, thick] (1.6,.4) -- (1.8,.6) -- (1.6,.8);
\node[probranch] {}
    child {
        node[probranch] {} 
           child {
        				node[probranch] {}    
        				edge from parent 
           					node[] {}
            				node[]  {}
            }     
        edge from parent 
            node[] {}
            node[]  {}
    }
    child {
        node[probranch] {}
        		child {
        				node[probranch] {}    
        				edge from parent 
           					node[] {}
            				node[]  {}
    				}    
        edge from parent 
            node[] {}
            node[]  {}
    }
    child {
        node[probranch, label=above: $1$] {}
        		child {
        				node[probranch] {}    
        				edge from parent 
           					node[] {}
            				node[]  {}
    				}
    				child {
        				node[probranch, label=above: $c$] {}    
        				edge from parent 
           					node[] {}
            				node[]  {}
    				}          
        edge from parent         
            node[] {}
            node[]  {}
    };
\end{tikzpicture}\\
\end{tabular}
\tikzstyle{level 1}=[level distance=9mm, sibling distance=6mm]
\tikzstyle{level 2}=[level distance=6mm, sibling distance=4mm]
\begin{tabular}{l r}
\begin{tikzpicture}[grow=right, sloped]
\draw  [-, rounded corners, thick] (1.6,.4) -- (1.8,.6) -- (1.6,.8);
\node[probranch] {}
    child {
        node[probranch] {} 
           child {
        				node[probranch] {}    
        				edge from parent 
           					node[] {}
            				node[]  {}
            }     
        edge from parent 
            node[] {}
            node[]  {}
    }
    child {
        node[probranch] {}
        		child {
        				node[probranch] {}    
        				edge from parent 
           					node[] {}
            				node[]  {}
    				}    
        edge from parent 
            node[] {}
            node[]  {}
    }
    child {
        node[probranch, label=above: $1$] {}
        		child {
        				node[probranch] {}    
        				edge from parent 
           					node[] {}
            				node[]  {}
    				}
    				child {
        				node[probranch, label=above: $d$] {}    
        				edge from parent 
           					node[] {}
            				node[]  {}
    				}          
        edge from parent         
            node[] {}
            node[]  {}
    };
\end{tikzpicture}
\end{tabular}\\
\end{center}
\begin{center}
\tikzstyle{level 1}=[level distance=9mm, sibling distance=7mm]
\tikzstyle{level 2}=[level distance=6mm, sibling distance=4mm]
\begin{tabular}{l r}
\begin{tikzpicture}[grow=right, sloped]
\node[probranch] {}
    child {
        node[probranch] {}
        		child {
        				node[probranch] {}    
        				edge from parent 
           					node[] {}
            				node[]  {}
    				}
    				child {
        				node[probranch] {}    
        				edge from parent 
           					node[] {}
            				node[]  {}
    				}     
        edge from parent 
            node[] {}
            node[]  {}
    }
    child {
        node[probranch, label=above: $1$] {}
        		child {
        				node[probranch] {}    
        				edge from parent 
           					node[] {}
            				node[]  {}
    				} 
    				child {
        				node[probranch, label=above: $e$] {}    
        				edge from parent 
           					node[] {}
            				node[]  {}
    				}          
        edge from parent         
            node[] {}
            node[]  {}
    };
\end{tikzpicture} 
\end{tabular}
\tikzstyle{level 1}=[level distance=9mm, sibling distance=7mm]
\tikzstyle{level 2}=[level distance=6mm, sibling distance=4mm]
\begin{tabular}{l r}
\begin{tikzpicture}[grow=right, sloped]
\draw  [-, rounded corners, thick] (1.6,.15) -- (1.8,.35) -- (1.6,.55);
\node[probranch] {}
    child {
        node[probranch] {}
        		child {
        				node[probranch] {}    
        				edge from parent 
           					node[] {}
            				node[]  {}
    				}
    				child {
        				node[probranch] {}    
        				edge from parent 
           					node[] {}
            				node[]  {}
    				}     
        edge from parent 
            node[] {}
            node[]  {}
    }
    child {
        node[probranch, label=above: $1$] {}
        		child {
        				node[probranch] {}    
        				edge from parent 
           					node[] {}
            				node[]  {}
    				} 
    				child {
        				node[probranch, label=above: $f$] {}    
        				edge from parent 
           					node[] {}
            				node[]  {}
    				}          
        edge from parent         
            node[] {}
            node[]  {}
    };
\end{tikzpicture} \\
\end{tabular}
\tikzstyle{level 1}=[level distance=9mm, sibling distance=7mm]
\tikzstyle{level 2}=[level distance=6mm, sibling distance=4mm]
\begin{tabular}{l r}
\begin{tikzpicture}[grow=right, sloped]
\draw  [-, rounded corners, thick] (1.6,-.15) -- (1.8,-.35) -- (1.6,-.55);
\draw  [-, rounded corners, thick] (1.6,.15) -- (1.8,.35) -- (1.6,.55);
\node[probranch] {}
    child {
        node[probranch] {}
        		child {
        				node[probranch] {}    
        				edge from parent 
           					node[] {}
            				node[]  {}
    				}
    				child {
        				node[probranch] {}    
        				edge from parent 
           					node[] {}
            				node[]  {}
    				}     
        edge from parent 
            node[] {}
            node[]  {}
    }
    child {
        node[probranch, label=above: $1$] {}
        		child {
        				node[probranch] {}    
        				edge from parent 
           					node[] {}
            				node[]  {}
    				} 
    				child {
        				node[probranch, label=above: $g$] {}    
        				edge from parent 
           					node[] {}
            				node[]  {}
    				}          
        edge from parent         
            node[] {}
            node[]  {}
    };
\end{tikzpicture}\\
\end{tabular}
\end{center}
$a={\left\{\frac{n}{2}\right\}}^{15?}_4$, $b={\left\{\frac{n}{2}\right\}}^{15?}_4$, $c={\left\{n\right\}}^{7}_2$, $d={\left\{n\right\}}^7_2$, $e={\left\{\frac{n}{2}\right\}}^{8}_4$, $f={\left\{n\right\}}^{4}_2$, $g={\left\{n\right\}}^{4}_2$ \\ 
\\
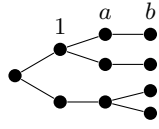
\begin{center}
\tikzstyle{level 1}=[level distance=6mm, sibling distance=7mm]
\tikzstyle{level 2}=[level distance=6mm, sibling distance=4mm]
\tikzstyle{level 3}=[level distance=6mm, sibling distance=3mm]
\begin{tabular}{l r}
\begin{tikzpicture}[grow=right, sloped]
\node[probranch] {}
    child {
        node[probranch] {}
        		child {
        				node[probranch] {}
        						child {
        								node[probranch] {}    
        						edge from parent 
           							node[] {}
            						node[]  {}
    								}
    								child {
        								node[probranch] {}    
        						edge from parent 
           							node[] {}
            						node[]  {}
    								}      
        				edge from parent 
           					node[] {}
            				node[]  {}
    				}   
        edge from parent 
            node[] {}
            node[]  {}
    }
    child {
        node[probranch, label=above: $1$] {}
        		child {
        				node[probranch] {}
        						child {
        								node[probranch] {}    
        						edge from parent 
           							node[] {}
            						node[]  {}
    								}    
        				edge from parent 
           					node[] {}
            				node[]  {}
    				} 
    				child {
        				node[probranch, label=above: $a$] {}
        						child {
        								node[probranch, label=above: $b$] {}    
        						edge from parent 
           							node[] {}
            						node[]  {}
    								}    
        				edge from parent 
           					node[] {}
            				node[]  {}
    				}          
        edge from parent         
            node[] {}
            node[]  {}
    };
\end{tikzpicture} \\
\end{tabular}
\end{center}
$\left\{(a,b)\right\}=\left\{(\frac{3}{2},x):x=\frac{4}{2},\frac{5}{2},...,\frac{14}{2}\right\}\cup\left\{(2,x):x=\frac{5}{2},\frac{6}{2},...,\frac{13}{2}\right\}\cup\left\{(\frac{5}{2},x):x=\frac{6}{2},\frac{7}{2},...,\frac{12}{2}\right\}\cup\left\{(3,x):x=\frac{7}{2},\frac{8}{2},...,\frac{11}{2}\right\}\cup\left\{(\frac{7}{2},x):x=\frac{8}{2},\frac{9}{2},\frac{10}{2}\right\}\cup\left\{(4,\frac{4}{9})\right\} $\\
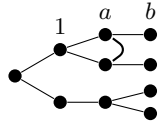
\begin{center}
\tikzstyle{level 1}=[level distance=6mm, sibling distance=7mm]
\tikzstyle{level 2}=[level distance=6mm, sibling distance=4mm]
\tikzstyle{level 3}=[level distance=6mm, sibling distance=3mm]
\begin{tabular}{l r}
\begin{tikzpicture}[grow=right, sloped]
\draw  [-, rounded corners, thick] (1.3,.2) -- (1.5,.35) -- (1.3,.5);
\node[probranch] {}
    child {
        node[probranch] {}
        		child {
        				node[probranch] {}
        						child {
        								node[probranch] {}    
        						edge from parent 
           							node[] {}
            						node[]  {}
    								}
    								child {
        								node[probranch] {}    
        						edge from parent 
           							node[] {}
            						node[]  {}
    								}      
        				edge from parent 
           					node[] {}
            				node[]  {}
    				}   
        edge from parent 
            node[] {}
            node[]  {}
    }
    child {
        node[probranch, label=above: $1$] {}
        		child {
        				node[probranch] {}
        						child {
        								node[probranch] {}    
        						edge from parent 
           							node[] {}
            						node[]  {}
    								}    
        				edge from parent 
           					node[] {}
            				node[]  {}
    				} 
    				child {
        				node[probranch, label=above: $a$] {}
        						child {
        								node[probranch, label=above: $b$] {}    
        						edge from parent 
           							node[] {}
            						node[]  {}
    								}    
        				edge from parent 
           					node[] {}
            				node[]  {}
    				}          
        edge from parent         
            node[] {}
            node[]  {}
    };
\end{tikzpicture} \\
\end{tabular}
\end{center}
$\left\{(a,b)\right\}=\left\{(2,x):x=\frac{5}{2},\frac{6}{2},...,\frac{13}{2}\right\}\cup\left\{(3,x):x=\frac{7}{2},\frac{8}{2},...,\frac{11}{2}\right\}\cup\left\{(4,\frac{4}{9})\right\} $
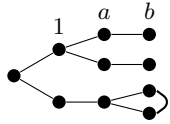
\begin{center}
\tikzstyle{level 1}=[level distance=6mm, sibling distance=7mm]
\tikzstyle{level 2}=[level distance=6mm, sibling distance=4mm]
\tikzstyle{level 3}=[level distance=6mm, sibling distance=3mm]
\begin{tabular}{l r}
\begin{tikzpicture}[grow=right, sloped]
\draw  [-, rounded corners, thick] (1.9,-.2) -- (2.1,-.35) -- (1.9,-.5);
\node[probranch] {}
    child {
        node[probranch] {}
        		child {
        				node[probranch] {}
        						child {
        								node[probranch] {}    
        						edge from parent 
           							node[] {}
            						node[]  {}
    								}
    								child {
        								node[probranch] {}    
        						edge from parent 
           							node[] {}
            						node[]  {}
    								}      
        				edge from parent 
           					node[] {}
            				node[]  {}
    				}   
        edge from parent 
            node[] {}
            node[]  {}
    }
    child {
        node[probranch, label=above: $1$] {}
        		child {
        				node[probranch] {}
        						child {
        								node[probranch] {}    
        						edge from parent 
           							node[] {}
            						node[]  {}
    								}    
        				edge from parent 
           					node[] {}
            				node[]  {}
    				} 
    				child {
        				node[probranch, label=above: $a$] {}
        						child {
        								node[probranch, label=above: $b$] {}    
        						edge from parent 
           							node[] {}
            						node[]  {}
    								}    
        				edge from parent 
           					node[] {}
            				node[]  {}
    				}          
        edge from parent         
            node[] {}
            node[]  {}
    };
\end{tikzpicture} \\
\end{tabular}

\end{center}
$\left\{(a,b)\right\}=\left\{(\frac{3}{2},x):x=2,3,...7\right\}\cup\left\{(2,x):x=3,4,5,6\right\}\cup\left\{(\frac{5}{2},x):x=3,4,5,6\right\}\cup\left\{(3,x):x=4,5\right\}\cup\left\{(\frac{7}{2},x):x=4,5\right\}$\\
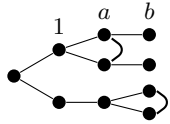
\begin{center}
\tikzstyle{level 1}=[level distance=6mm, sibling distance=7mm]
\tikzstyle{level 2}=[level distance=6mm, sibling distance=4mm]
\tikzstyle{level 3}=[level distance=6mm, sibling distance=3mm]
\begin{tabular}{l r}
\begin{tikzpicture}[grow=right, sloped]
\draw  [-, rounded corners, thick] (1.3,.2) -- (1.5,.35) -- (1.3,.5);
\draw  [-, rounded corners, thick] (1.9,-.2) -- (2.1,-.35) -- (1.9,-.5);
\node[probranch] {}
    child {
        node[probranch] {}
        		child {
        				node[probranch] {}
        						child {
        								node[probranch] {}    
        						edge from parent 
           							node[] {}
            						node[]  {}
    								}
    								child {
        								node[probranch] {}    
        						edge from parent 
           							node[] {}
            						node[]  {}
    								}      
        				edge from parent 
           					node[] {}
            				node[]  {}
    				}   
        edge from parent 
            node[] {}
            node[]  {}
    }
    child {
        node[probranch, label=above: $1$] {}
        		child {
        				node[probranch] {}
        						child {
        								node[probranch] {}    
        						edge from parent 
           							node[] {}
            						node[]  {}
    								}    
        				edge from parent 
           					node[] {}
            				node[]  {}
    				} 
    				child {
        				node[probranch, label=above: $a$] {}
        						child {
        								node[probranch, label=above: $b$] {}    
        						edge from parent 
           							node[] {}
            						node[]  {}
    								}    
        				edge from parent 
           					node[] {}
            				node[]  {}
    				}          
        edge from parent         
            node[] {}
            node[]  {}
    };
\end{tikzpicture} \\
\end{tabular}
\end{center}
$\left\{(a,b)\right\}=\left\{(2,x):x=3,4,5,6\right\}\cup\left\{(3,x):x=4,5\right\} $ \\
\\
Tangent Cones: \\ $y(x-y)(x+y)(x+ay)$ where $a\neq0,1,-1$, \\ $y(x^2+y^2)(x+ay)$ where $a\neq0$ \\
\\
Newton Polygons 38, 39, 40: \\
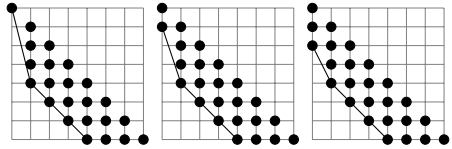
\begin{center}
\begin{tikzpicture}[scale=.25]
\draw[help lines] (0,0) grid (7,7);
\draw(0,7)--(1,3);
\draw(1,3)--(4,0);
\draw[fill] (0,7) circle [radius= 0.25];
\draw[fill] (1,6) circle [radius= 0.25];
\draw[fill] (2,5) circle [radius= 0.25];
\draw[fill] (3,4) circle [radius= 0.25];
\draw[fill] (4,3) circle [radius= 0.25];
\draw[fill] (5,2) circle [radius= 0.25];
\draw[fill] (6,1) circle [radius= 0.25];
\draw[fill] (7,0) circle [radius= 0.25];
\draw[fill] (6,0) circle [radius= 0.25];
\draw[fill](5,1) circle [radius = 0.25];
\draw[fill](4,2) circle [radius= 0.25];
\draw[fill](3,3) circle [radius = 0.25];
\draw[fill](2,4) circle[radius = 0.25];
\draw[fill](5,0) circle [radius = 0.25];
\draw[fill](4,1) circle [radius = 0.25];
\draw[fill](3,2) circle [radius = 0.25];
\draw[fill](4,0) circle [radius = 0.25];
\draw[fill](1,5) circle [radius = 0.25];
\draw[fill](2,3) circle [radius = 0.25];
\draw[fill](3,1) circle [radius = 0.25];
\draw[fill](2,2) circle [radius = 0.25];
\draw[fill](1,4) circle [radius = 0.25];
\draw[fill](1,3) circle [radius = 0.25];
\end{tikzpicture}
\begin{tikzpicture}[scale=.25]
\draw[help lines] (0,0) grid (7,7);
\draw(0,6)--(1,3);
\draw(1,3)--(4,0);
\draw[fill] (0,7) circle [radius= 0.25];
\draw[fill] (1,6) circle [radius= 0.25];
\draw[fill] (2,5) circle [radius= 0.25];
\draw[fill] (3,4) circle [radius= 0.25];
\draw[fill] (4,3) circle [radius= 0.25];
\draw[fill] (5,2) circle [radius= 0.25];
\draw[fill] (6,1) circle [radius= 0.25];
\draw[fill] (7,0) circle [radius= 0.25];
\draw[fill] (6,0) circle [radius= 0.25];
\draw[fill](5,1) circle [radius = 0.25];
\draw[fill](4,2) circle [radius= 0.25];
\draw[fill](3,3) circle [radius = 0.25];
\draw[fill](2,4) circle[radius = 0.25];
\draw[fill](5,0) circle [radius = 0.25];
\draw[fill](4,1) circle [radius = 0.25];
\draw[fill](3,2) circle [radius = 0.25];
\draw[fill](4,0) circle [radius = 0.25];
\draw[fill](1,5) circle [radius = 0.25];
\draw[fill](2,3) circle [radius = 0.25];
\draw[fill](3,1) circle [radius = 0.25];
\draw[fill](2,2) circle [radius = 0.25];
\draw[fill](1,4) circle [radius = 0.25];
\draw[fill](0,6) circle [radius = 0.25];
\draw[fill](1,3) circle [radius = 0.25];
\end{tikzpicture}
\begin{tikzpicture}[scale=.25]
\draw[help lines] (0,0) grid (7,7);
\draw(0,5)--(1,3);
\draw(1,3)--(4,0);
\draw[fill] (0,7) circle [radius= 0.25];
\draw[fill] (1,6) circle [radius= 0.25];
\draw[fill] (2,5) circle [radius= 0.25];
\draw[fill] (3,4) circle [radius= 0.25];
\draw[fill] (4,3) circle [radius= 0.25];
\draw[fill] (5,2) circle [radius= 0.25];
\draw[fill] (6,1) circle [radius= 0.25];
\draw[fill] (7,0) circle [radius= 0.25];
\draw[fill] (6,0) circle [radius= 0.25];
\draw[fill](5,1) circle [radius = 0.25];
\draw[fill](4,2) circle [radius= 0.25];
\draw[fill](3,3) circle [radius = 0.25];
\draw[fill](2,4) circle[radius = 0.25];
\draw[fill](5,0) circle [radius = 0.25];
\draw[fill](4,1) circle [radius = 0.25];
\draw[fill](3,2) circle [radius = 0.25];
\draw[fill](4,0) circle [radius = 0.25];
\draw[fill](1,5) circle [radius = 0.25];
\draw[fill](2,3) circle [radius = 0.25];
\draw[fill](3,1) circle [radius = 0.25];
\draw[fill](2,2) circle [radius = 0.25];
\draw[fill](1,4) circle [radius = 0.25];
\draw[fill](0,6) circle [radius = 0.25];
\draw[fill](0,5) circle [radius = 0.25];
\draw[fill](1,3) circle [radius = 0.25];
\end{tikzpicture}
\end{center}
Families: \\ $y(x-y)(x+y)(x+ay)+bx^4y+cx^3y^2+dx^2y^3+exy^4+fy^5+gx^5y+hx^4y^2+jx^3y^3+kx^2y^4+lxy^5+my^6+nx^7+px^6y+qx^5y^2+rx^4y^3+sx^3y^4+tx^2y^5+uxy^6+vy^7=0$, where $a\neq0,1,-1$\\ \\ $y(x^2+y^2)(x+ay)+cx^4y+dx^3y^2+ex^2y^3+fxy^4+gy^5+hx^5y+jx^4y^2+kx^3y^3+lx^2y^4+mxy^5+ny^6+px^7+qx^6y+rx^5y^2+sx^4y^3+tx^3y^4+ux^2y^5+vxy^6+wy^7=0$,  where $a\neq0$ \\ \\ $y(x-y)(x+y)(x+ay)+cx^4y+dx^3y^2+ex^2y^3+fxy^4+gy^5+hx^6+jx^5y+kx^4y^2+lx^3y^3+mx^2y^4+nxy^5+py^6+qx^7+rx^6y+sx^5y^2+tx^4y^3+ux^3y^4+vx^2y^5+wxy^6+zy^7=0$, where $a\neq0,1,-1$ \\ \\ $y(x^2+y^2)(x+ay)+cx^4y+dx^3y^2+ex^2y^3+fxy^4+gy^5+hx^6+jx^5y+kx^4y^2+lx^3y^3+mx^2y^4+nxy^5+py^6+qx^7+rx^6y+sx^5y^2+tx^4y^3+ux^3y^4+vx^2y^5+wxy^6+zy^7=0$,  where $a\neq0$ \\ \\ $y(x-y)(x+y)(x+ay)+cx^5+dx^4y+ex^3y^2+fx^2y^3+gxy^4+hy^5+jx^6+kx^5y+lx^4y^2+mx^3y^3+nx^2y^4+pxy^5+qy^6+rx^7+sx^6y+tx^5y^2+ux^4y^3+vx^3y^4+wx^2y^5+zxy^6+Ay^7=0$, where $a\neq0,1,-1$ \\ \\ $y(x^2+y^2)(x+ay)+cx^5+dx^4y+ex^3y^2+fx^2y^3+gxy^4+hy^5+jx^6+kx^5y+lx^4y^2+mx^3y^3+nx^2y^4+pxy^5+qy^6+rx^7+sx^6y+tx^5y^2+ux^4y^3+vx^3y^4+wx^2y^5+zxy^6+{\alpha}y^7=0$  where $a\neq0$ \\ \\
Diagrams: 
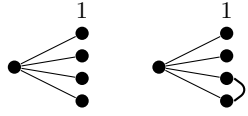
\begin{center}
\tikzstyle{level 1}=[level distance=9mm, sibling distance=3mm]
\tikzstyle{level 2}=[level distance=6mm, sibling distance=6mm]
\begin{tabular}{l r}
\begin{tikzpicture}[grow=right, sloped]
\node[probranch] {}
    child {
        node[probranch] {}        
        edge from parent 
            node[] {}
            node[]  {}
    }
    child {
        node[probranch] {}        
        edge from parent 
            node[] {}
            node[]  {}
    }
    child {
        node[probranch] {}        
        edge from parent 
            node[] {}
            node[]  {}
    }
    child {
        node[probranch, label=above: $1$] {}        
        edge from parent         
            node[] {}
            node[]  {}
    };
\end{tikzpicture}\\
\end{tabular}
\tikzstyle{level 1}=[level distance=9mm, sibling distance=3mm]
\tikzstyle{level 2}=[level distance=6mm, sibling distance=6mm]
\begin{tabular}{l r}
\begin{tikzpicture}[grow=right, sloped]
\draw  [-, rounded corners, thick] (1,-.15) -- (1.2,-.3) -- (1,-.45);
\node[probranch] {}
    child {
        node[probranch] {}        
        edge from parent 
            node[] {}
            node[]  {}
    }
    child {
        node[probranch] {}        
        edge from parent 
            node[] {}
            node[]  {}
    }
    child {
        node[probranch] {}        
        edge from parent 
            node[] {}
            node[]  {}
    }
    child {
        node[probranch, label=above: $1$] {}        
        edge from parent         
            node[] {}
            node[]  {}
    };
\end{tikzpicture}\\
\end{tabular}
\end{center}
Only two new diagrams. All other families transform to Newton polygons 24 or 25 or 26 or 27 or 28 or 29 or 30 or 31. \\ 
\\
Tangent Cones: \\ $(y^2+x^2)(y^2+ax^2)$ where $a>0$ and $a\neq1$, \\ $(y^2+x^2)^2$ \\
\\
Newton Polygon 41: \\ 
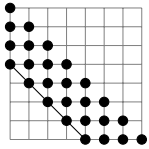
\begin{center}
\begin{tikzpicture}[scale=.25]
\draw[help lines] (0,0) grid (7,7);
\draw(0,4)--(4,0);
\draw[fill] (0,7) circle [radius= 0.25];
\draw[fill] (1,6) circle [radius= 0.25];
\draw[fill] (2,5) circle [radius= 0.25];
\draw[fill] (3,4) circle [radius= 0.25];
\draw[fill] (4,3) circle [radius= 0.25];
\draw[fill] (5,2) circle [radius= 0.25];
\draw[fill] (6,1) circle [radius= 0.25];
\draw[fill] (7,0) circle [radius= 0.25];
\draw[fill] (6,0) circle [radius= 0.25];
\draw[fill](5,1) circle [radius = 0.25];
\draw[fill](4,2) circle [radius= 0.25];
\draw[fill](3,3) circle [radius = 0.25];
\draw[fill](2,4) circle[radius = 0.25];
\draw[fill](5,0) circle [radius = 0.25];
\draw[fill](4,1) circle [radius = 0.25];
\draw[fill](3,2) circle [radius = 0.25];
\draw[fill](4,0) circle [radius = 0.25];
\draw[fill](1,5) circle [radius = 0.25];
\draw[fill](2,3) circle [radius = 0.25];
\draw[fill](3,1) circle [radius = 0.25];
\draw[fill](2,2) circle [radius = 0.25];
\draw[fill](1,4) circle [radius = 0.25];
\draw[fill](0,6) circle [radius = 0.25];
\draw[fill](0,5) circle [radius = 0.25];
\draw[fill](1,3) circle [radius = 0.25];
\draw[fill](0,4) circle [radius = 0.25];
\end{tikzpicture}
\end{center}
Families: \\ $(y^2+x^2)(y^2+ax^2)+bx^5+cx^4y+dx^3y^2+ex^2y^3+fxy^4+gy^5+hx^6+jx^5y+kx^4y^2+lx^3y^3+mx^2y^4+nxy^5+py^6+q^7+rx^6y+sx^5y^2+tx^4y^3+ux^3y^4+vx^2y^5+wxy^6+zy^7=0$ where $a>0$ and $a\neq1$, \\ \\ $(y^2+x^2)^2+ax^5+bx^4y+cx^3y^2+dx^2y^3+exy^4+fy^5+gx^6+hx^5y+jx^4y^2+kx^3y^3+lx^2y^4+mxy^5+ny^6+px^7+qx^6y+rx^5y^2+sx^4y^3+tx^3y^4+ux^2y^5+vxy^6+wy^7=0$ \\ \\
Diagrams: 
\begin{center}
\tikzstyle{level 1}=[level distance=9mm, sibling distance=3mm]
\tikzstyle{level 2}=[level distance=6mm, sibling distance=6mm]
\begin{center}
Next, we show the corresponding diagram types\\
\end{center}

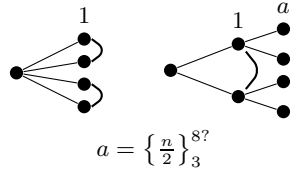
\begin{tabular}{l r}
\begin{tikzpicture}[grow=right, sloped]
\draw  [-, rounded corners, thick] (1,.15) -- (1.2,.3) -- (1,.45);
\draw  [-, rounded corners, thick] (1,-.15) -- (1.2,-.3) -- (1,-.45);
\node[probranch] {}
    child {
        node[probranch] {}        
        edge from parent 
            node[] {}
            node[]  {}
    }
    child {
        node[probranch] {}        
        edge from parent 
            node[] {}
            node[]  {}
    }
    child {
        node[probranch] {}        
        edge from parent 
            node[] {}
            node[]  {}
    }
    child {
        node[probranch, label=above: $1$] {}        
        edge from parent         
            node[] {}
            node[]  {}
    };
\end{tikzpicture}\\
\end{tabular}
\tikzstyle{level 1}=[level distance=9mm, sibling distance=7mm]
\tikzstyle{level 2}=[level distance=6mm, sibling distance=4mm]
\begin{tabular}{l r}
\begin{tikzpicture}[grow=right, sloped]
\draw  [-, rounded corners, thick] (1,.25) -- (1.2,0) -- (1,-.25);
\node[probranch] {}
    child {
        node[probranch] {}
        		child {
        				node[probranch] {}    
        				edge from parent 
           					node[] {}
            				node[]  {}
    				}
    				child {
        				node[probranch] {}    
        				edge from parent 
           					node[] {}
            				node[]  {}
    				}     
        edge from parent 
            node[] {}
            node[]  {}
    }
    child {
        node[probranch, label=above: $1$] {}
        		child {
        				node[probranch] {}    
        				edge from parent 
           					node[] {}
            				node[]  {}
    				} 
    				child {
        				node[probranch, label=above: $a$] {}    
        				edge from parent 
           					node[] {}
            				node[]  {}
    				}          
        edge from parent         
            node[] {}
            node[]  {}
    };
\end{tikzpicture}\\
\end{tabular}

$a={\left\{\frac{n}{2}\right\}}^{8?}_3$
\end{center}

\section{Open Problems}

At first glance it may seem surprising that the lower the multiplicity, the harder it is to classify the singular points for curves of fixed degree.  But on second thought, the equation for a family of curves having a singular point of lower multiplicity has more terms, so it stands to reason that the calculation is more difficult.  The system of diagrams that we have developed enables one to visualize, in some way, all possible types of singular points of algebraic curves.  In this section, we wish to describe some open problems in the classification of singular points. \\ \\

\subsection{Double points}  

For curves of degree 7 determine all values of $n$ and $m$ for which there exists a curve having a singular point with diagram \\

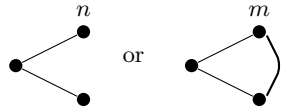
\begin{center}
\tikzstyle{level 1}=[level distance=9mm, sibling distance=9mm]
\tikzstyle{level 2}=[level distance=6mm, sibling distance=6mm]
\begin{tabular}{l r}
\begin{tikzpicture}[grow=right, sloped]
\node[probranch] {}
    child 
    {
        node[probranch]{}        
        edge from parent 
            node[] {}
            node[]  {}
    }
    child 
    {
        node[probranch, label=above: $n$] {}
        edge from parent         
        	  node[] {}
        	  node[]  {}
    };
\end{tikzpicture}  \\
\end{tabular}
or
\tikzstyle{level 1}=[level distance=9mm, sibling distance=9mm]
\tikzstyle{level 2}=[level distance=6mm, sibling distance=6mm]
\begin{tabular}{l r}
\begin{tikzpicture}[grow=right, sloped]
\draw  [-, rounded corners, thick] (1,.4) -- (1.2,0) -- (1,-.4);
\node[probranch] {}
    child 
    {
        node[probranch]{}        
        edge from parent 
            node[] {}
            node[]  {}
    }
    child 
    {
        node[probranch, label=above: $m$] {}
        edge from parent         
        	  node[] {}
        	  node[]  {}
    };
\end{tikzpicture}\\
\end{tabular}
\end{center}

(For complex curves, one only considers the diagram without braces.)  In general, for curves of degree $d$, determine all values of $n$ and $m$ for which there exists a curve having a singular point with diagram \\

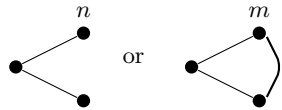
\begin{center}
\tikzstyle{level 1}=[level distance=9mm, sibling distance=9mm]
\tikzstyle{level 2}=[level distance=6mm, sibling distance=6mm]
\begin{tabular}{l r}
\begin{tikzpicture}[grow=right, sloped]
\node[probranch] {}
    child 
    {
        node[probranch]{}        
        edge from parent 
            node[] {}
            node[]  {}
    }
    child 
    {
        node[probranch, label=above: $n$] {}
        edge from parent         
        	  node[] {}
        	  node[]  {}
    };
\end{tikzpicture}  \\
\end{tabular}
or
\tikzstyle{level 1}=[level distance=9mm, sibling distance=9mm]
\tikzstyle{level 2}=[level distance=6mm, sibling distance=6mm]
\begin{tabular}{l r}
\begin{tikzpicture}[grow=right, sloped]
\draw  [-, rounded corners, thick] (1,.4) -- (1.2,0) -- (1,-.4);
\node[probranch] {}
    child 
    {
        node[probranch]{}        
        edge from parent 
            node[] {}
            node[]  {}
    }
    child 
    {
        node[probranch, label=above: $m$] {}
        edge from parent         
        	  node[] {}
        	  node[]  {}
    };
\end{tikzpicture}\\
\end{tabular}
\end{center}

(In Arnol'd's notation, singular points of type $A_n$ or $A_m^*$.)\\

\subsection{Triple points}

For curves of degree 7, determine all values of $n$ and $m$ for which there exists a curve having a singular point with diagram \\

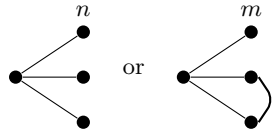
\begin{center}
\tikzstyle{level 1}=[level distance=9mm, sibling distance=6mm]
\begin{tabular}{l r}
\begin{tikzpicture}[grow=right, sloped]
\node[probranch] {}
    child {
        node[probranch] {}        
        edge from parent 
            node[] {}
            node[]  {}
    }
    child {
        node[probranch] {}        
        edge from parent 
            node[] {}
            node[]  {}
    }
    child {
        node[probranch, label=above: $n$] {}        
        edge from parent         
            node[] {}
            node[]  {}
    };
\end{tikzpicture}\\
\end{tabular}
or
\tikzstyle{level 1}=[level distance=9mm, sibling distance=6mm]
\begin{tabular}{l r}
\begin{tikzpicture}[grow=right, sloped]
\draw  [-, rounded corners, thick] (1,0) -- (1.2,-.3) -- (1,-.6);
\node[probranch] {}
    child {
        node[probranch] {}        
        edge from parent 
            node[] {}
            node[]  {}
    }
    child {
        node[probranch] {}        
        edge from parent 
            node[] {}
            node[]  {}
    }
    child {
        node[probranch, label=above: $m$] {}        
        edge from parent         
            node[] {}
            node[]  {}
    };
\end{tikzpicture}\\
\end{tabular}

\end{center}
and also determine all values of $q,r,s,$ and $t$

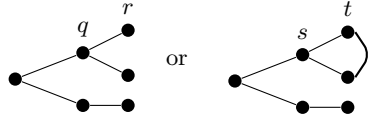
\begin{center}
\tikzstyle{level 1}=[level distance=9mm, sibling distance=7mm]
\tikzstyle{level 2}=[level distance=6mm, sibling distance=6mm]
\begin{tabular}{l r}
\begin{tikzpicture}[grow=right, sloped]
\node[probranch] {}
    child {
        node[probranch] {}
        		child {
        				node[probranch] {}    
        				edge from parent 
           					node[] {}
            				node[]  {}
    				}    
        edge from parent 
            node[] {}
            node[]  {}
    }
    child {
        node[probranch, label=above: $q$] {}
        		child {
        				node[probranch] {}    
        				edge from parent 
           					node[] {}
            				node[]  {}
    				}
    				child {
        				node[probranch, label=above: $r$] {}    
        				edge from parent 
           					node[] {}
            				node[]  {}
    				}          
        edge from parent         
            node[] {}
            node[]  {}
    };
\end{tikzpicture} \\
\end{tabular}
or
\tikzstyle{level 1}=[level distance=9mm, sibling distance=7mm]
\tikzstyle{level 2}=[level distance=6mm, sibling distance=6mm]
\begin{tabular}{l r}
\begin{tikzpicture}[grow=right, sloped]
\draw  [-, rounded corners, thick] (1.6,0.05) -- (1.8,.35) -- (1.6,.65);
\node[probranch] {}
    child {
        node[probranch] {}
        		child {
        				node[probranch] {}    
        				edge from parent 
           					node[] {}
            				node[]  {}
    				}    
        edge from parent 
            node[] {}
            node[]  {}
    }
    child {
        node[probranch, label=above: $s$] {}
        		child {
        				node[probranch] {}    
        				edge from parent 
           					node[] {}
            				node[]  {}
    				}
    				child {
        				node[probranch, label=above: $t$] {}    
        				edge from parent 
           					node[] {}
            				node[]  {}
    				}          
        edge from parent         
            node[] {}
            node[]  {}
    };
\end{tikzpicture}\\
\end{tabular}
\end{center}
such that there exists a curve having a singular point with one of the diagrams above.  You could ask the same questions for algebraic curves of degree $d$.

\subsection{Quadruple points}

For curves of degree 7, we found that for $a= \frac{7}{4},\frac{8}{4},...,\frac{17}{4}$ and $b=2,3,4$ there exist curves having diagrams of type 

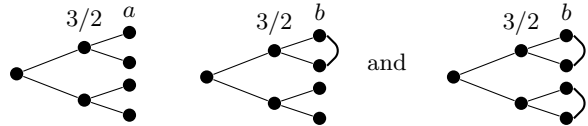
\begin{center}
\tikzstyle{level 1}=[level distance=9mm, sibling distance=7mm]
\tikzstyle{level 2}=[level distance=6mm, sibling distance=4mm]
\begin{tabular}{l r}
\begin{tikzpicture}[grow=right, sloped]
\node[probranch] {}
    child {
        node[probranch] {}
        		child {
        				node[probranch] {}    
        				edge from parent 
           					node[] {}
            				node[]  {}
    				}
    				child {
        				node[probranch] {}    
        				edge from parent 
           					node[] {}
            				node[]  {}
    				}     
        edge from parent 
            node[] {}
            node[]  {}
    }
    child {
        node[probranch, label=above: $3/2$] {}
        		child {
        				node[probranch] {}    
        				edge from parent 
           					node[] {}
            				node[]  {}
    				} 
    				child {
        				node[probranch, label=above: $a$] {}    
        				edge from parent 
           					node[] {}
            				node[]  {}
    				}          
        edge from parent         
            node[] {}
            node[]  {}
    };
\end{tikzpicture} \\
\end{tabular}
\tikzstyle{level 1}=[level distance=9mm, sibling distance=7mm]
\tikzstyle{level 2}=[level distance=6mm, sibling distance=4mm]
\begin{tabular}{l r}
\begin{tikzpicture}[grow=right, sloped]
\draw  [-, rounded corners, thick] (1.6,.15) -- (1.8,.35) -- (1.6,.55);
\node[probranch] {}
    child {
        node[probranch] {}
        		child {
        				node[probranch] {}    
        				edge from parent 
           					node[] {}
            				node[]  {}
    				}
    				child {
        				node[probranch] {}    
        				edge from parent 
           					node[] {}
            				node[]  {}
    				}     
        edge from parent 
            node[] {}
            node[]  {}
    }
    child {
        node[probranch, label=above: $3/2$] {}
        		child {
        				node[probranch] {}    
        				edge from parent 
           					node[] {}
            				node[]  {}
    				} 
    				child {
        				node[probranch, label=above: $b$] {}    
        				edge from parent 
           					node[] {}
            				node[]  {}
    				}          
        edge from parent         
            node[] {}
            node[]  {}
    };
\end{tikzpicture} \\
\end{tabular}
and
\tikzstyle{level 1}=[level distance=9mm, sibling distance=7mm]
\tikzstyle{level 2}=[level distance=6mm, sibling distance=4mm]
\begin{tabular}{l r}
\begin{tikzpicture}[grow=right, sloped]
\draw  [-, rounded corners, thick] (1.6,-.15) -- (1.8,-.35) -- (1.6,-.55);
\draw  [-, rounded corners, thick] (1.6,.15) -- (1.8,.35) -- (1.6,.55);
\node[probranch] {}
    child {
        node[probranch] {}
        		child {
        				node[probranch] {}    
        				edge from parent 
           					node[] {}
            				node[]  {}
    				}
    				child {
        				node[probranch] {}    
        				edge from parent 
           					node[] {}
            				node[]  {}
    				}     
        edge from parent 
            node[] {}
            node[]  {}
    }
    child {
        node[probranch, label=above: $3/2$] {}
        		child {
        				node[probranch] {}    
        				edge from parent 
           					node[] {}
            				node[]  {}
    				} 
    				child {
        				node[probranch, label=above: $b$] {}    
        				edge from parent 
           					node[] {}
            				node[]  {}
    				}          
        edge from parent         
            node[] {}
            node[]  {}
    };
\end{tikzpicture} \\
\end{tabular}
\end{center}
Because of the size of the computer calculation, we were not able to show that these were the only values of $a$ and $b$ for 7th degree curves.  Also, for curves of degree 7, we found that for $a=\frac{5}{2},\frac{6}{2},...,\frac{11}{2}$ and $b=3,4,5$ there exist curves having diagrams \\

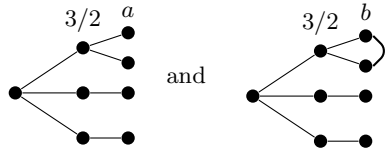
\begin{center}
\tikzstyle{level 1}=[level distance=9mm, sibling distance=6mm]
\tikzstyle{level 2}=[level distance=6mm, sibling distance=4mm]
\begin{tabular}{l r}
\begin{tikzpicture}[grow=right, sloped]
\node[probranch] {}
    child {
        node[probranch] {} 
           child {
        				node[probranch] {}    
        				edge from parent 
           					node[] {}
            				node[]  {}
            }     
        edge from parent 
            node[] {}
            node[]  {}
    }
    child {
        node[probranch] {}
        		child {
        				node[probranch] {}    
        				edge from parent 
           					node[] {}
            				node[]  {}
    				}    
        edge from parent 
            node[] {}
            node[]  {}
    }
    child {
        node[probranch, label=above: $3/2$] {}
        		child {
        				node[probranch] {}    
        				edge from parent 
           					node[] {}
            				node[]  {}
    				}
    				child {
        				node[probranch, label=above: $a$] {}    
        				edge from parent 
           					node[] {}
            				node[]  {}
    				}          
        edge from parent         
            node[] {}
            node[]  {}
    };
\end{tikzpicture} \\
\end{tabular}
and
\tikzstyle{level 1}=[level distance=9mm, sibling distance=6mm]
\tikzstyle{level 2}=[level distance=6mm, sibling distance=4mm]
\begin{tabular}{l r}
\begin{tikzpicture}[grow=right, sloped]
\draw  [-, rounded corners, thick] (1.6,0.4) -- (1.8,.6) -- (1.6,.8);
\node[probranch] {}
    child {
        node[probranch] {} 
           child {
        				node[probranch] {}    
        				edge from parent 
           					node[] {}
            				node[]  {}
            }     
        edge from parent 
            node[] {}
            node[]  {}
    }
    child {
        node[probranch] {}
        		child {
        				node[probranch] {}    
        				edge from parent 
           					node[] {}
            				node[]  {}
    				}    
        edge from parent 
            node[] {}
            node[]  {}
    }
    child {
        node[probranch, label=above: $3/2$] {}
        		child {
        				node[probranch] {}    
        				edge from parent 
           					node[] {}
            				node[]  {}
    				}
    				child {
        				node[probranch, label=above: $b$] {}    
        				edge from parent 
           					node[] {}
            				node[]  {}
    				}          
        edge from parent         
            node[] {}
            node[]  {}
    };
\end{tikzpicture} \\
\end{tabular}
\end{center}
Again, because of the size of the computer calculation, we were not able to show that these are the only values of $a$ and $b$ for 7th degree curves.  Similarly, for 7th degree curves we found the values $d=\frac{5}{2},\frac{6}{2},...,\frac{14}{2}$ and $e=3,4,5,6,7$ for diagrams of the type \\

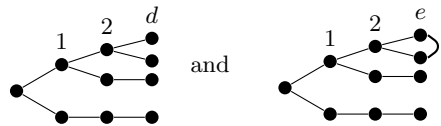
\begin{center}
\tikzstyle{level 1}=[level distance=6mm, sibling distance=7mm]
\tikzstyle{level 2}=[level distance=6mm, sibling distance=4mm]
\tikzstyle{level 3}=[level distance=6mm, sibling distance=3mm]
\begin{tabular}{l r}
\begin{tikzpicture}[grow=right, sloped]
\node[probranch] {}
    child {
        node[probranch] {}
        		child {
        				node[probranch] {}
        						child {
        								node[probranch] {}    
        						edge from parent 
           							node[] {}
            						node[]  {}
    								}   
        				edge from parent 
           					node[] {}
            				node[]  {}
    				}   
        edge from parent 
            node[] {}
            node[]  {}
    }
    child {
        node[probranch, label=above: $1$] {}
        		child {
        				node[probranch] {}
        						child {
        								node[probranch] {}    
        						edge from parent 
           							node[] {}
            						node[]  {}
    								}    
        				edge from parent 
           					node[] {}
            				node[]  {}
    				} 
    				child {
        				node[probranch, label=above: $2$] {}
        						child {
        								node[probranch] {}    
        						edge from parent 
           							node[] {}
            						node[]  {}
    								} 
        						child {
        								node[probranch, label=above: $d$] {}    
        						edge from parent 
           							node[] {}
            						node[]  {}
    								}    
        				edge from parent 
           					node[] {}
            				node[]  {}
    				}          
        edge from parent         
            node[] {}
            node[]  {}
    };
\end{tikzpicture} \\
\end{tabular}
and
\tikzstyle{level 1}=[level distance=6mm, sibling distance=7mm]
\tikzstyle{level 2}=[level distance=6mm, sibling distance=4mm]
\tikzstyle{level 3}=[level distance=6mm, sibling distance=3mm]
\begin{tabular}{l r}
\begin{tikzpicture}[grow=right, sloped]
\draw  [-, rounded corners, thick] (1.9,0.4) -- (2.1,.55) -- (1.9,.7);
\node[probranch] {}
    child {
        node[probranch] {}
        		child {
        				node[probranch] {}
        						child {
        								node[probranch] {}    
        						edge from parent 
           							node[] {}
            						node[]  {}
    								}   
        				edge from parent 
           					node[] {}
            				node[]  {}
    				}   
        edge from parent 
            node[] {}
            node[]  {}
    }
    child {
        node[probranch, label=above: $1$] {}
        		child {
        				node[probranch] {}
        						child {
        								node[probranch] {}    
        						edge from parent 
           							node[] {}
            						node[]  {}
    								}    
        				edge from parent 
           					node[] {}
            				node[]  {}
    				} 
    				child {
        				node[probranch, label=above: $2$] {}
        						child {
        								node[probranch] {}    
        						edge from parent 
           							node[] {}
            						node[]  {}
    								} 
        						child {
        								node[probranch, label=above: $e$] {}    
        						edge from parent 
           							node[] {}
            						node[]  {}
    								}    
        				edge from parent 
           					node[] {}
            				node[]  {}
    				}          
        edge from parent         
            node[] {}
            node[]  {}
    };
\end{tikzpicture} \\
\end{tabular}
\end{center}
but we were not able to show that these are the only values. \\

By now the reader should be aware that there are an unlimited number of open problems of the type:  for curves of degree $d$ find all values at the top of the columns of a given diagram for which there exists a curve of degree $d$ having a singular point with that diagram. \\

\renewcommand{\baselinestretch}{1}

\textsc{David A. Weinberg: Department of Mathematics and Statistics, Texas Tech University, Lubbock, Texas 79409-1042}

{\em E-mail address}: david.weinberg@ttu.edu\\

\textsc{Nicholas J. Willis: Department of Mathematics, George Fox University, Newberg, Oregon 97132}

{\em E-mail address}: nwillis@georgefox.edu \\

\end{document}